\documentclass[11pt,a4paper]{article}

\usepackage{array}
\newcolumntype{"}{@{\hskip\tabcolsep\vrule width 1.5pt\hskip\tabcolsep}}

\usepackage{graphics}
\usepackage{epsfig}
\usepackage{subfigure}
\usepackage{amssymb}
\usepackage{amsthm}
\usepackage{mathrsfs}
\usepackage{array}
\usepackage{hhline}
\usepackage{longtable}
\usepackage{amsmath}
\usepackage{float}
\usepackage{picinpar}
\usepackage{cite}
\usepackage{enumerate}
\usepackage{makecell}
\usepackage{multirow}
\newcommand{\PreserveBackslash}[1]{\let\temp=\\#1\let\\=\temp}
\let\PBS=\PreserveBackslash
\usepackage[mathlines]{lineno}
\newcommand*\patchAmsMathEnvironmentForLineno[1]{%
\expandafter\let\csname old#1\expandafter\endcsname\csname #1\endcsname  \expandafter\let\csname oldend#1\expandafter\endcsname\csname end#1\endcsname  \renewenvironment{#1}%
{\linenomath\csname old#1\endcsname}%
{\csname oldend#1\endcsname\endlinenomath}}%
\newcommand*\patchBothAmsMathEnvironmentsForLineno[1]{%
\patchAmsMathEnvironmentForLineno{#1}%
\patchAmsMathEnvironmentForLineno{#1*}}%
\AtBeginDocument{%
\patchBothAmsMathEnvironmentsForLineno{equation}%
\patchBothAmsMathEnvironmentsForLineno{align}%
\patchBothAmsMathEnvironmentsForLineno{flalign}%
\patchBothAmsMathEnvironmentsForLineno{alignat}%
\patchBothAmsMathEnvironmentsForLineno{gather}%
\patchBothAmsMathEnvironmentsForLineno{multline}%
}

\newcommand{\mb}[1]{\text{\boldmath ${#1}$}}
\def\Z{\mathbb{Z}}

\def\d{\delta}
\def\D{\Delta}

\usepackage[dvipsnames]{xcolor}
\definecolor{indiagreen}{rgb}{0.07, 0.53, 0.03}
\definecolor{vividviolet}{rgb}{0.62, 0.0, 1.0}
\def\D{\Delta}
\def\rsq{\hspace*{\fill}$\blacksquare$\medskip}

\newtheorem{theorem}{Theorem}[section]
\newtheorem{lemma}[theorem]{Lemma}
\newtheorem{corollary}[theorem]{Corollary}
\newtheorem{problem}{Problem}[section]

\newtheorem{rem}{Remark}[section]
\numberwithin{equation}{section}

\newtheoremstyle{example}
  {10pt}          
  {10pt}  
  {\rm}  
  {}
  {\bf}  
  {: }    
  { }    
  {}     
\theoremstyle{example}
\newtheorem{example}{Example}[section]
\def\Z{\mathbb{Z}}

\newtheorem{defi}{Definition}[section]

\newenvironment{definition}{\begin{defi}\rm }{\end{defi}}

\newcommand{\T}[1]{\fontsize{#1}{#1}\selectfont}

\usepackage{graphicx}
\graphicspath{%
    {converted_graphics/}
    {/}
}

\def\ms{\medskip}

\def\nt{\noindent}

\def\Shiu2{\color{vividviolet}}
\begin{document}

\begin{center}
{\mathversion{bold}\Large \bf Graphs With Minimal Strength}

\bigskip
{\large Zhen-Bin Gao$^a$, Gee-Choon Lau$^{b,}$\footnote{Corresponding author.}}, Wai-Chee Shiu{$^{c,d}$}, \\

\emph{{$^a$}College of Mathematical Sciences, Harbin Engineering University,}\\
\emph{Harbin 150001, P. R. China}\\
\emph{gaozhenbin@aliyun.com}\\

\medskip
\emph{{$^d$}Faculty of Computer \& Mathematical Sciences,}\\
\emph{Universiti Teknologi MARA (Segamat Campus),}\\
\emph{85000, Johor, Malaysia.}\\
\emph{geeclau@yahoo.com}\\

\medskip
\emph{{$^c$}College of Global Talents, Beijing Institute of Technology,\\ Zhuhai, China.}\\
\emph{{$^d$}Department of Mathematics, The Chinese University of Hong Kong,\\ Shatin, Hong Kong.}\\
\emph{wcshiu@associate.hkbu.edu.hk}
\end{center}
%

\begin{abstract}
For any graph $G$ of order $p$, a bijection $f: V(G)\to [1,p]$ is called a numbering of the graph $G$ of order $p$.  The strength $str_f(G)$ of a numbering $f: V(G)\to [1,p]$ of $G$ is defined by $str_f(G) = \max\{f(u)+f(v)\; |\; uv\in E(G)\},$ and the strength $str(G)$ of a graph $G$ itself is $str(G) = \min\{str_f(G)\;|\; f \mbox{ is a numbering of } G\}.$ A numbering $f$ is called a strength labeling of $G$ if $str_f(G)=str(G)$. In this paper, we obtained a sufficient condition for a graph to have $str(G)=|V(G)|+\d(G)$.  Consequently, many questions raised in [Bounds for the strength of graphs, {\it Aust. J. Combin.} {\bf72(3)},  (2018) 492--508] and [On the strength of some trees, {\it AKCE Int. J. Graphs Comb.} (Online 2019) doi.org/10.1016/j.akcej.2019.06.002] are solved. Moreover, we showed that every graph $G$ either has $str(G)=|V(G)|+\d(G)$ or is a proper subgraph of a graph $H$ that has $str(H) = |V(H)| + \d(H)$ with $\d(H)=\d(G)$. Further, new good lower bounds of $str(G)$ are also obtained. Using these, we determined the strength of 2-regular graphs and obtained new lower bounds of $str(Q_n)$ for various $n$, where $Q_n$ is the $n$-regular hypercube.

\ms

\end{abstract}


\section{Introduction}

We only consider simple and loopless $(p,q)$-graph $G=(V,E)$ with order $|V(G)|=p$ and size $|E(G)|=q$. If $uv\in E(G)$, we say $v$ is a neighbor of $u$. The degree of a vertex $v$ in a graph $G$ is the number of neighbors of $v$ in $G$, denoted $\deg_G(v)$. The {\it minimum degree} (and {\it maximum degree}) of $G$ is the {\it minimum} (and {\it maximum}) degree among the vertices of $G$, denoted $\d(G)$ (and $\D(G)$). A vertex of degree 0 is called an isolated vertex and a vertex of degree 1 is called a pendant vertex, with its incident edge is called a pendant edge. The set of all neighbors of $u$ is denoted $N_G(u)$. For $S\subset V(G)$, let $N_G(S)$ be the set of all neighbors of the vertices in $S$. We shall drop the subscript $G$ if no ambiguity.  For $a<b$, the set of integers from $a$ to $b$ is denoted $[a,b]$. All notation not defined in the paper are referred to~\cite{Bondy}.

\ms\nt The notion of strength of a graph $G$ was introduced in~\cite{Ichishima+MB+Oshima18} as a generalization of super magic strength~\cite{Avadayappan+Jeyanthi+Vasuki01} that is effectively defined only for super edge-magic graphs~\cite{Enomoto+Llado+Nakamigawa98} (also called strong vertex-graceful~\cite{Shiu+Wong10} and strongly indexable~\cite{Acharya+Hegde91}), to any nonempty graphs as follows.

\ms\nt A bijection $f: V(G)\to [1,p]$ is called a {\it numbering} of the graph $G$ of order $p$. 

\begin{definition} The {\it strength} $str_f(G)$ of a numbering $f: V(G)\to [1,p]$ of $G$ is defined by \[str_f(G) = \max\{f(u)+f(v)\; |\; uv\in E(G)\},\] and the {\it strength $str(G)$} of a graph $G$ itself is \[str(G) = \min\{str_f(G)\;|\; f \mbox{ is a numbering of } G\}.\]
A numbering $f$ is called a {\it strength labeling} of $G$ if $str_f(G)=str(G)$.
\end{definition}

\nt Several lower and upper bounds for $str(G)$ were obtained in \cite{Ichishima+MB+Oshima18}.

\begin{lemma}\label{lem-subgraph} If $H$ is a subgraph of a graph $G$, then $$str(H)\le str(G).$$
\end{lemma}

\begin{lemma}\label{lem-d>0} For every graph $G$ of order $p$ with $\d(G)\ge 1$, $$str(G)\ge p+\d(G).$$
\end{lemma}

\begin{lemma}\label{lem-connect} For every graph $G$ of order $p$ with $\d(G)\ge 1$,
\[str(G)\ge p+\kappa'(G)\ge p+\kappa(G),\] where $\kappa(G)$ and  $\kappa'(G)$ are the connectivity and the edge connectivity of $G$, respectively.

\end{lemma}

\begin{lemma}\label{lem-maxdeg} For every nonempty graph $G$ of order $p$, \[\D(G)+2 \le str(G)\le 2p-1.\]  \end{lemma}


\nt Let $G+H$ be the disjoint union of $G$ and $H$ with $V(G+H) = V(G)\cup V(H)$ and $E(G+H)=E(G)\cup E(H)$. The disjoint union of $m$ copies of $G$ is denoted $mG$. We first extend Lemma~\ref{lem-d>0} to graphs with isolated vertices.

\begin{lemma}\label{lem-isolate} Let $G$ be a graph with $\d(G)\ge 1$. If $m\ge 1$, then $$str(G+mK_1)=str(G).$$  \end{lemma}

\begin{proof}
Let $f$ be a strength labeling of the graph $G$. We extend $f$ to a numbering of $G+mK_1$ by assigning all $m$ isolated vertices by labels in $[p+1,p+m]$. Clearly, $str_f(G+mK_1)=str_f(G)=str(G)$. Hence $str(G+mK_1)\le str(G)$. Combining with Lemma~\ref{lem-subgraph} we have the lemma.
\end{proof}

\nt Thus, from now on, we only consider graphs without isolated vertices. In \cite{Ichishima+MB+Oshima18,Ichishima+MB+Oshima19}, the authors show that $str(G)=p+\d(G)$ if $G$ is a path, cycle, complete graph, complete bipartite graph, ladder graph, prism graph, M\"{o}bius ladder, book graph or $K_{m,n}\times K_2$ each of which has order $p$. Moreover, if $str(G)=p+\d(G)$ for a graph $G$ of order $p$ with $\d(G)\ge 1$, then $str(G\odot nK_1) = (n+1)p+1$, where $G\odot H$ is the corona product of $G$ and $H$. The following problems are posed.

\begin{problem}\label{pbm-gb} Find good bounds for the strength of a graph. \end{problem}

\begin{problem}\label{pbm-minstr} Find sufficient conditions for a graph $G$ of order $p$ with $\d(G)\ge 1$ to ensure $str(G)=p+\d(G)$. \end{problem}

\begin{problem}\label{pbm-lobster} For every lobster $T$, determine the exact value of $str(T)$. \end{problem}

\begin{problem}\label{pbm-Qn} For every integer $n\ge 3$, determine the strength of $Q_n$, the $n$-dimensional hypercube. \end{problem}

\section{Sufficient Condition}


\ms\nt Let $\mathcal G_1=G_1=G$ be a graph of order $p$ with $p-2\ge \d(G)=\d_1\ge 1$. Suppose $\mathcal G_i$, $i\ge 1$, is not $m_iK_1$ for $m_i\ge 1$ nor $m_iK_1 + K_r$ with $m_i\ge 0$ for some $r\ge 2$.  We may denote $\mathcal G_i$ by $m_iK_1+G_i$, where $m_i\ge 0$, and $\d_i=\d(G_i)\ge 1$. Let $\mathcal G_{i+1}$ be a graph obtained from $G_i$ by deleting the $m_i\ge 1$ isolated vertices that exist in $\mathcal G_i$, and a vertex of degree $\d_i$ together with all its neighbors in $G_i$. Continue the procedure until $\mathcal G_s$ is either $m_sK_1$ with $m_s\ge 1$ or $m_sK_1 + K_r$ with $m_s\ge 0$ for some $r,s\ge 2$. This sequence $\{\mathcal G_i\}_{i= 1}^s$ of subgraphs is called a {\it $\d$-sequence} of $G$. When $\mathcal G_s=m_sK_1$, we let $\d_s=0$ by convention. Let $\widetilde{y}_j(G)=m_j+1-\d_j$ for $1\le j\le s$.

\begin{example}\label{ex-example} Following are two examples to illustrate the above construction. The black vertex is the chosen vertex that will be deleted at each stage.

\begin{enumerate}[1.]

\item \hfill{}\\
$\begin{array}{ccccccccc}
\epsfig{file=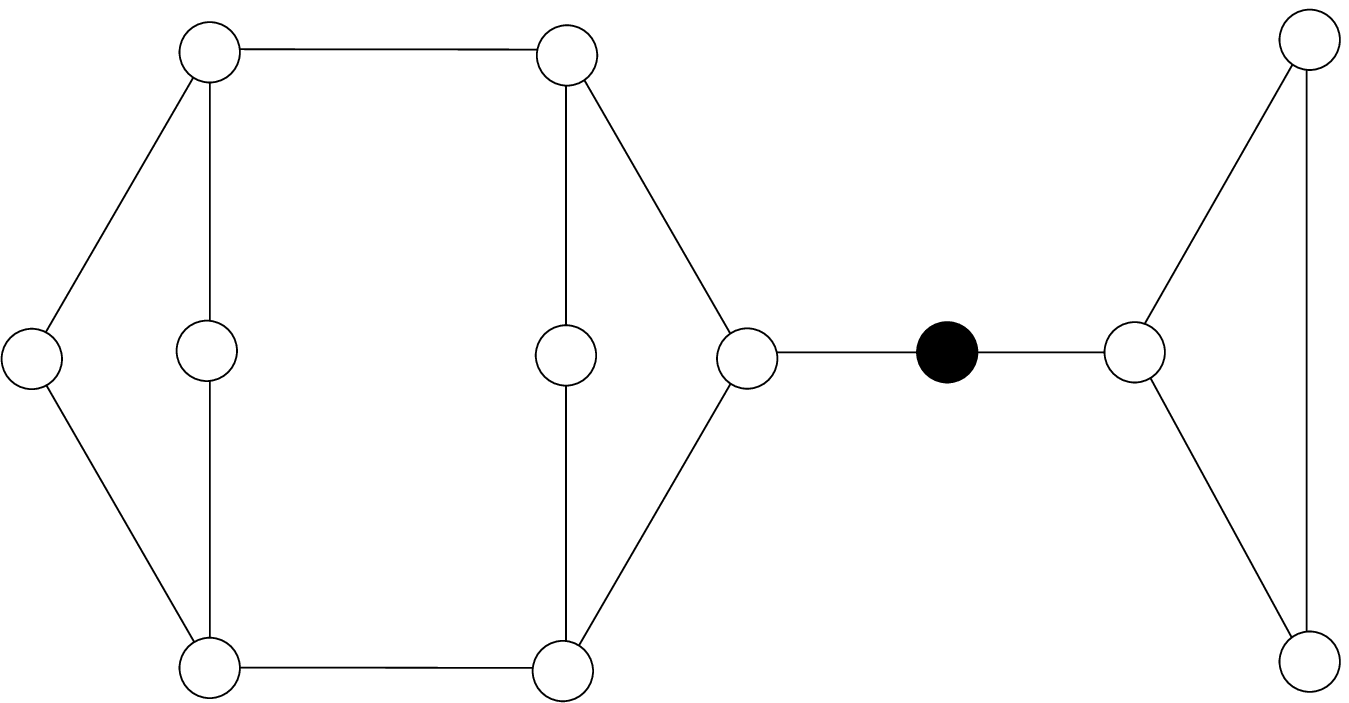, height=1.5cm} & \raisebox{6mm}[0mm][0mm]{$\rightarrow$} & \epsfig{file=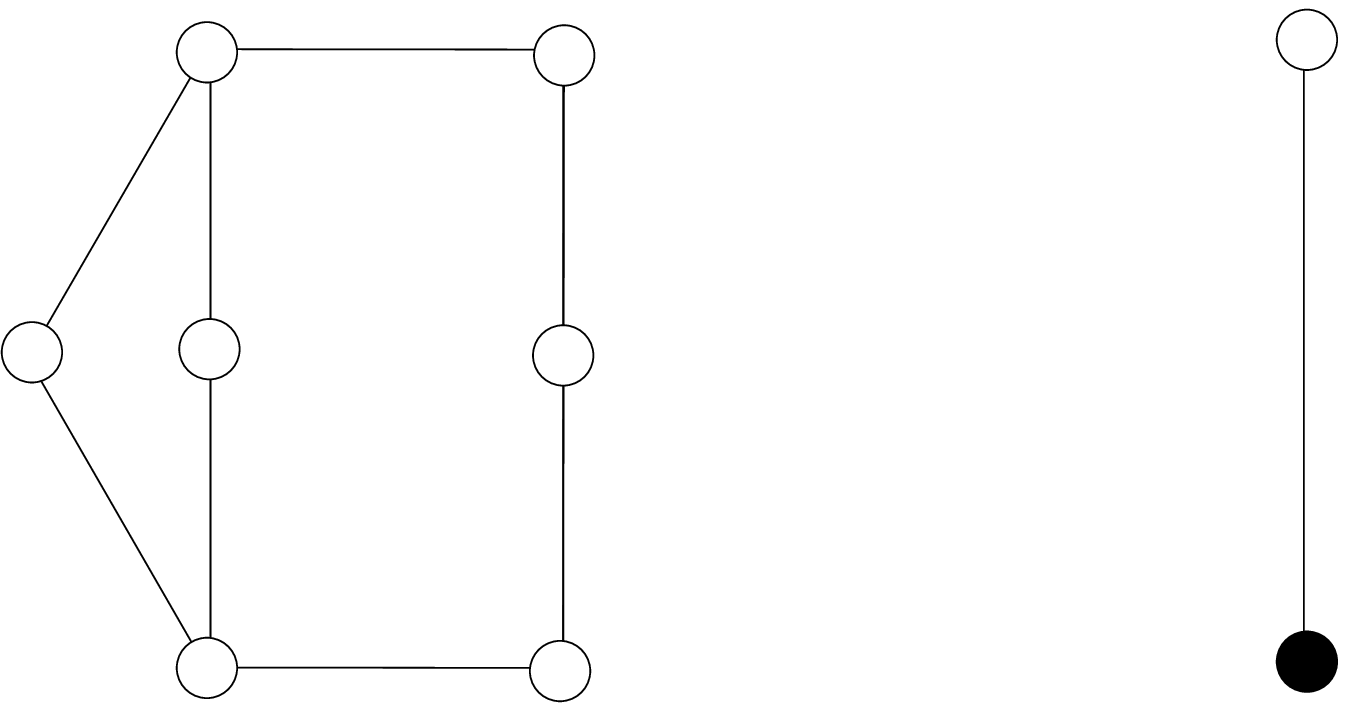, height=1.5cm}& \raisebox{6mm}[0mm][0mm]{$\rightarrow$} &\epsfig{file=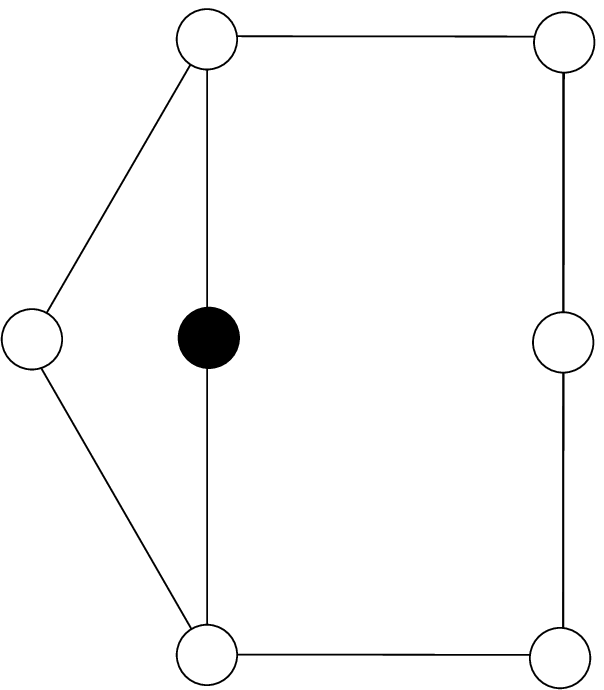, height=1.5cm} & \raisebox{6mm}[0mm][0mm]{$\rightarrow$} &\epsfig{file=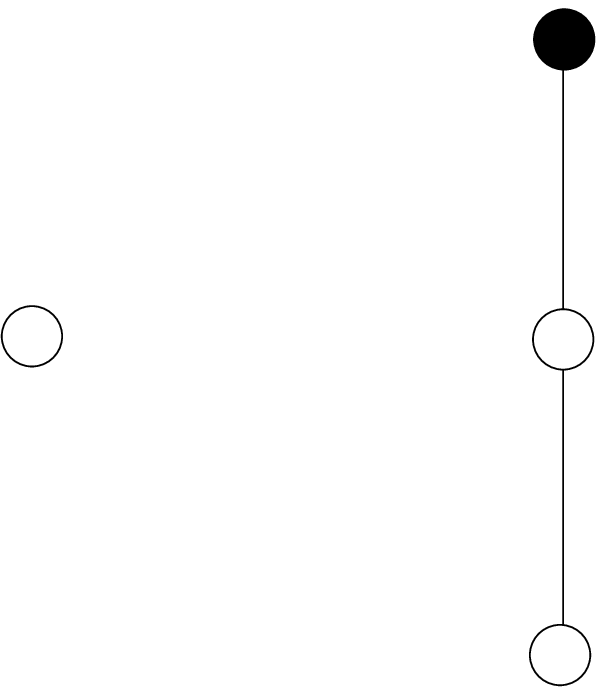, height=1.5cm} & \raisebox{6mm}[0mm][0mm]{$\rightarrow$} & \raisebox{6mm}[0mm][0mm]{\epsfig{file=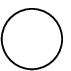, height=1.5mm}}\\
\mathcal G_1=G_1 & & \mathcal G_2=G_2 & & \mathcal G_3=G_3 & & \mathcal G_4=K_1+G_4 & & \mathcal G_5=K_1\\
\d_1=2 & & \d_2=1, m_2=0 & & \d_3=2, m_3=0 & & \d_4=1, m_4=1 & & \d_5=0, m_5=1
\end{array}$
\nt Now $(m_2+1-\d_2)+(m_3+1-\d_3)=(0)+(-1)=-1$. So, this $\d$-sequence of the graph $G_1$ does not satisfy the condition of \eqref{eq-cond} mentioned below.

\item \hfill{}\\
 $\begin{array}{ccccccc}
\epsfig{file=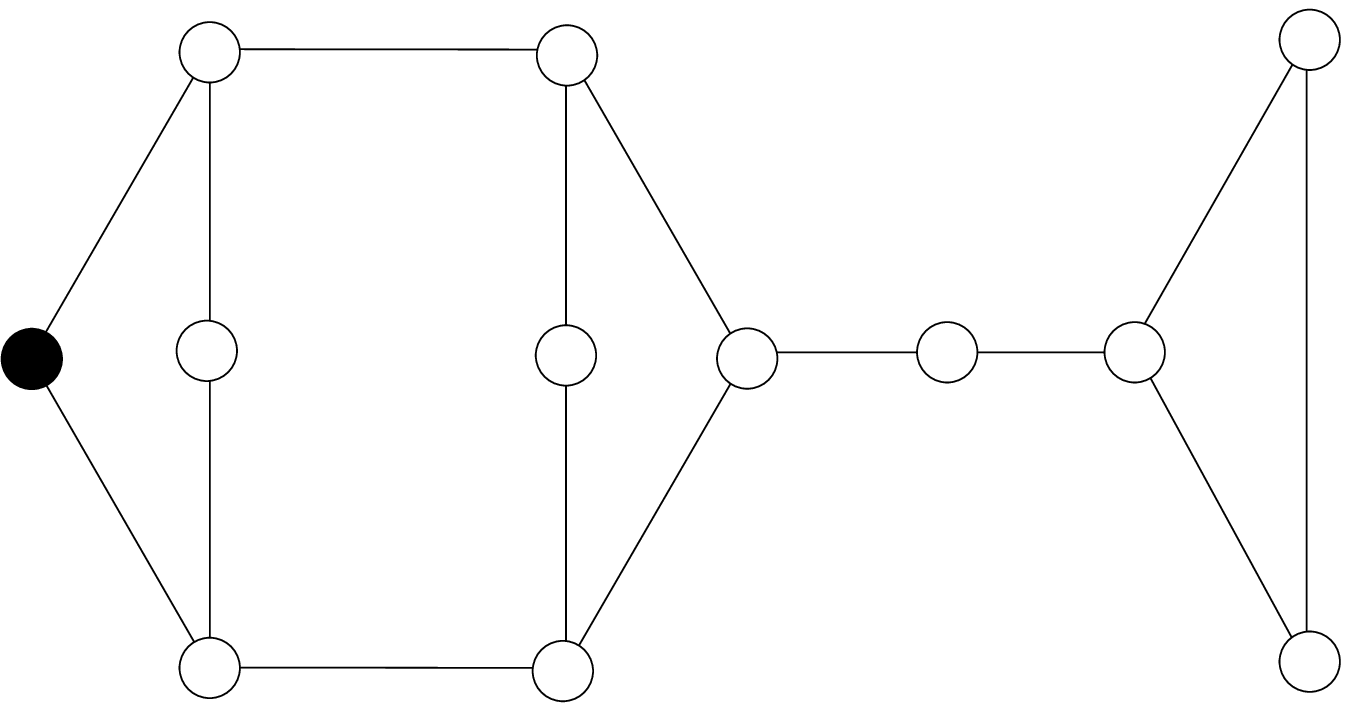, height=1.5cm} & \raisebox{6mm}[0mm][0mm]{$\rightarrow$} & \epsfig{file=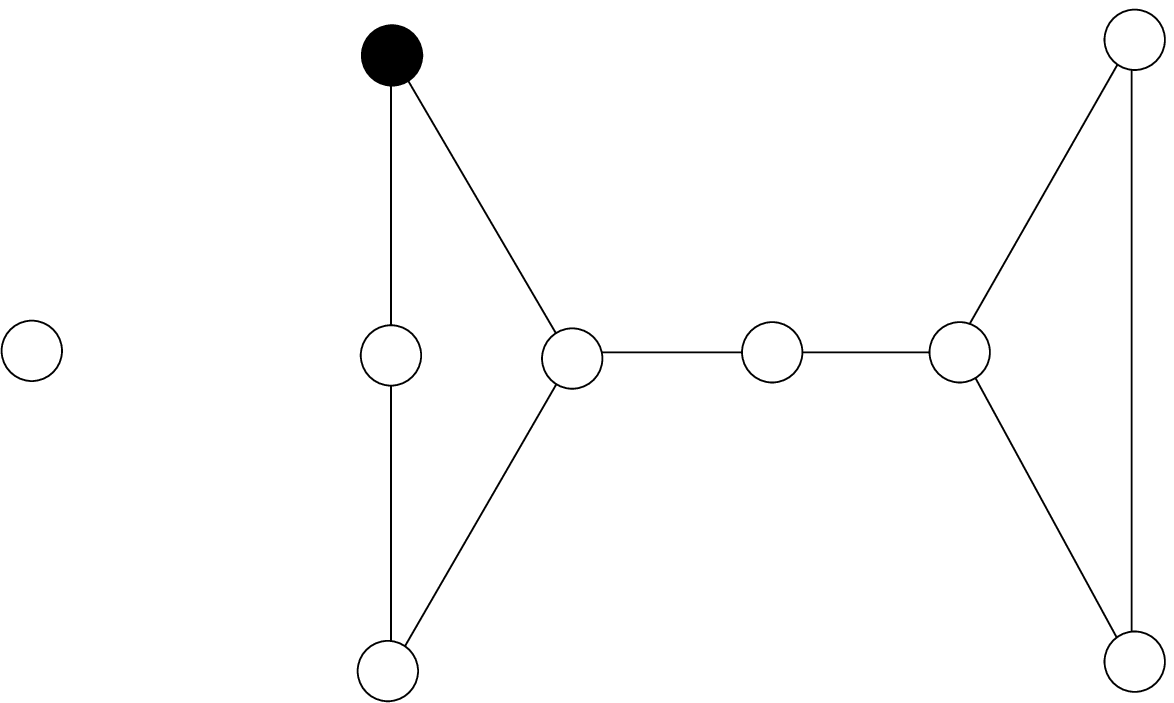, height=1.5cm}& \raisebox{6mm}[0mm][0mm]{$\rightarrow$} &\epsfig{file=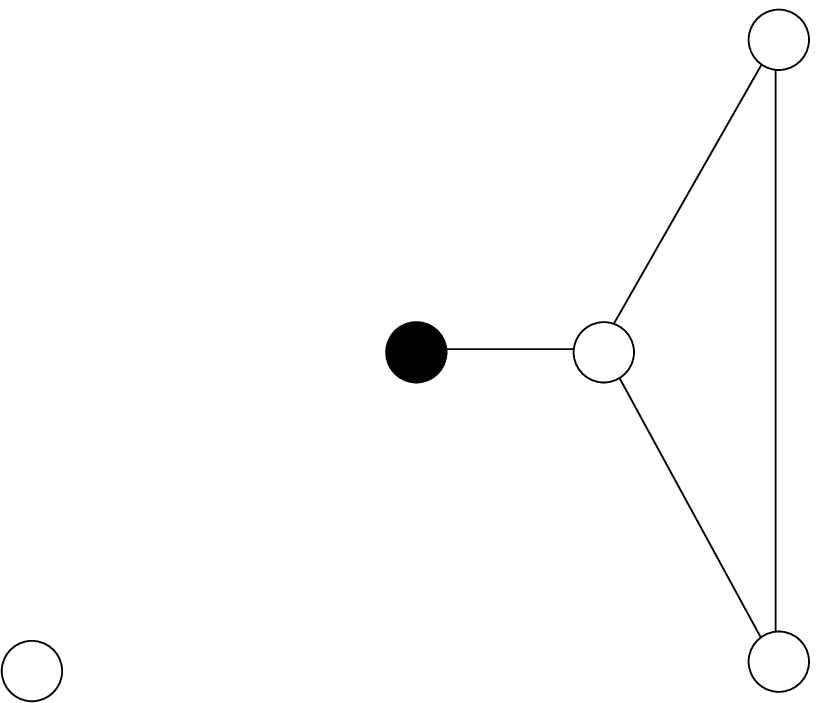, height=1.5cm} & \raisebox{6mm}[0mm][0mm]{$\rightarrow$} &\epsfig{file=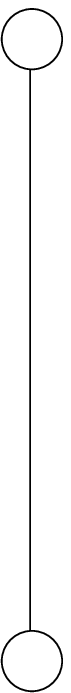, height=1.5cm} \\
\mathcal G_1=G_1 & & \mathcal G_2=K_1+G_2 & & \mathcal G_3=K_1+G_3 & & \mathcal G_4=K_2\\
\d_1=2 & & \d_2=2, m_2=1 & & \d_3=1, m_3=1 & & \d_4=1, m_4=0
\end{array}$

\nt This $\d$-sequence of the graph $G_1$ satisfies the condition of \eqref{eq-cond}. \rsq
\end{enumerate}

\end{example}

\begin{example}\label{ex-counterexample} Consider the following graph $G$:\\
\centerline{\epsfig{file=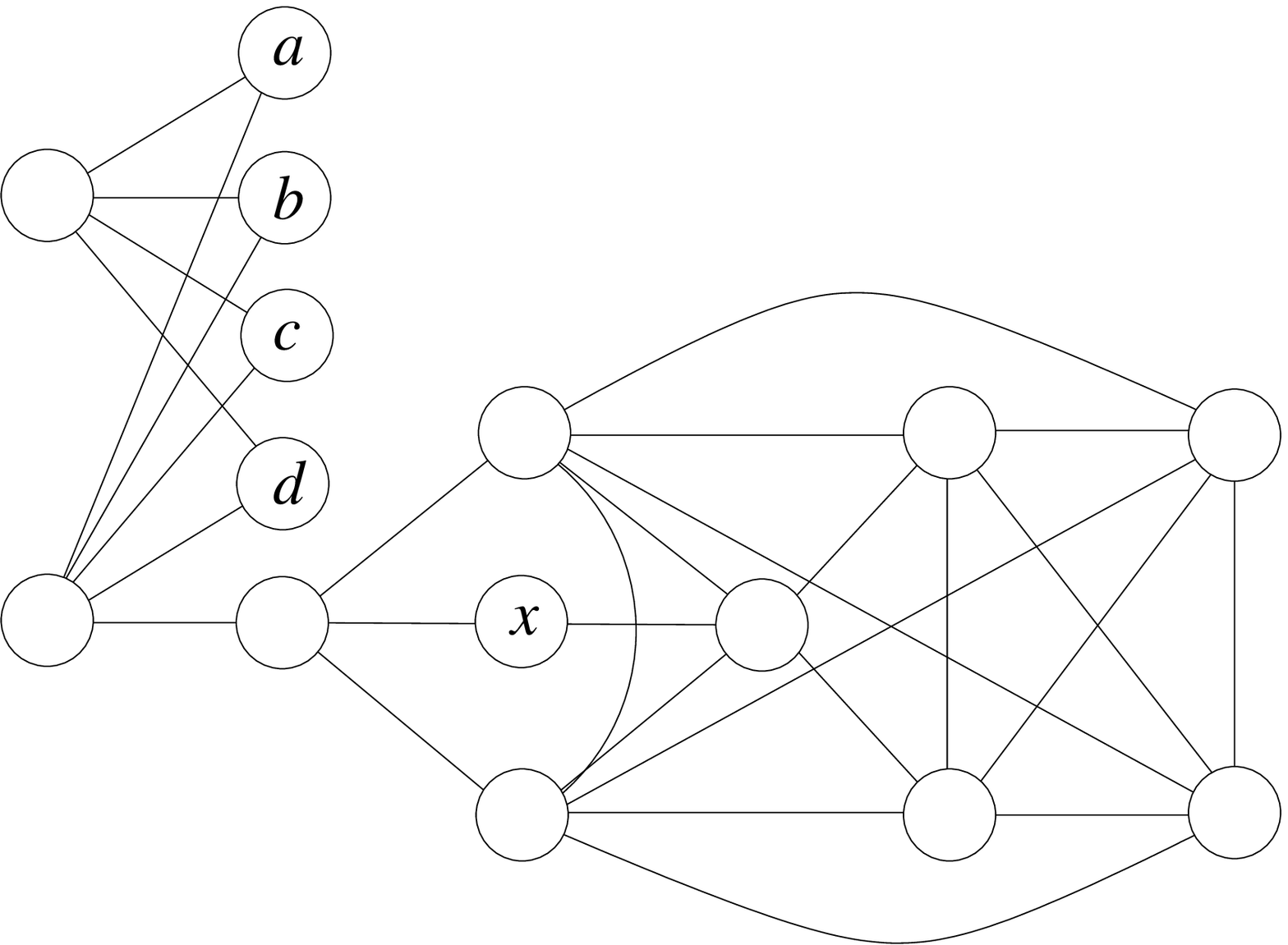, width=4cm}}

\nt There are 5 vertices of degree $2=\d(G)$.
\begin{enumerate}[(1).]
\item Suppose we choose $a$ as the first vertex (similarly can choosse $b$, $c$, or $d$). So we have

$\begin{array}{ccccccc}
\epsfig{file=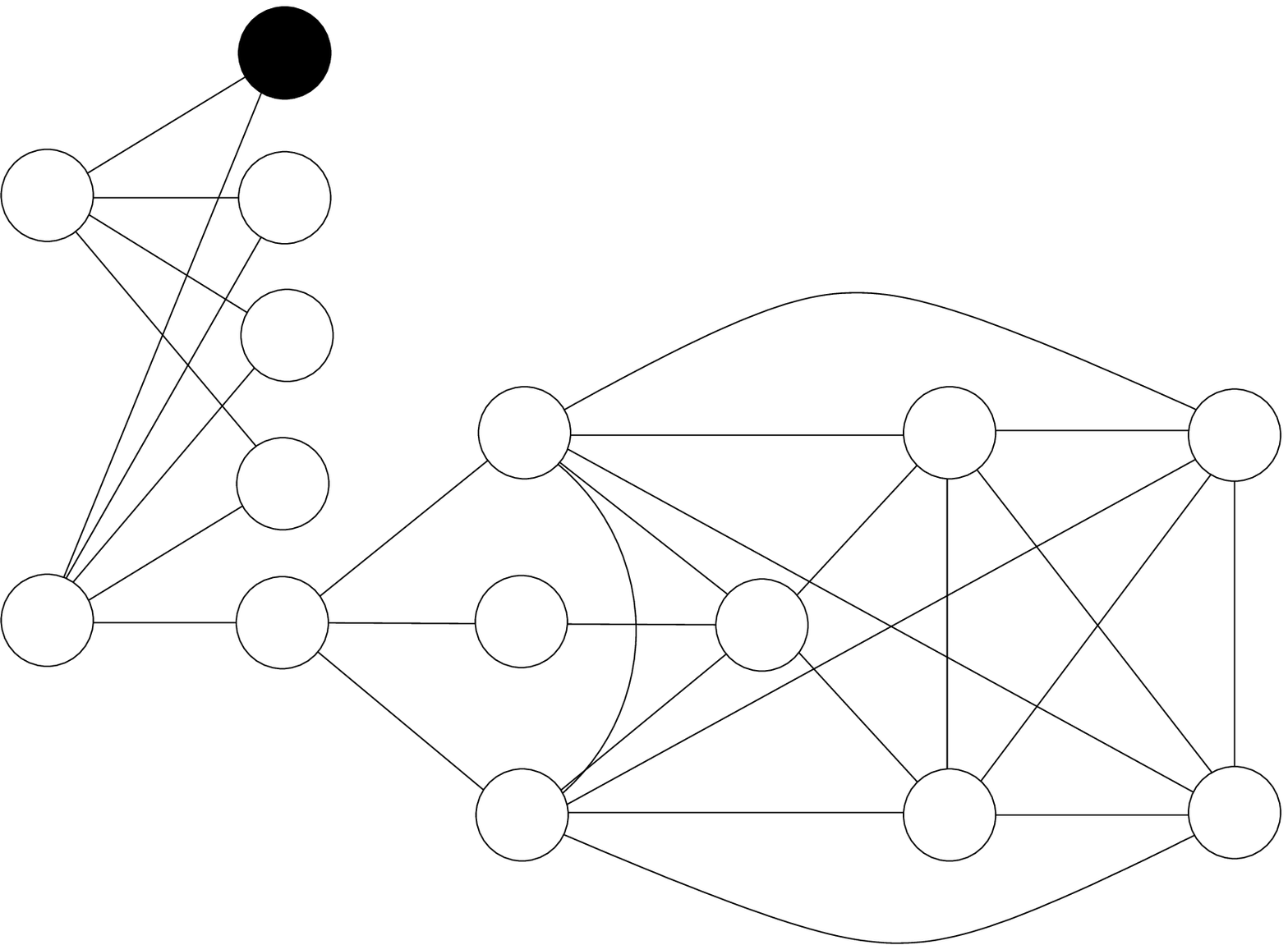, width=3cm} &\raisebox{10mm}[0mm][0mm]{$\rightarrow$}&
\epsfig{file=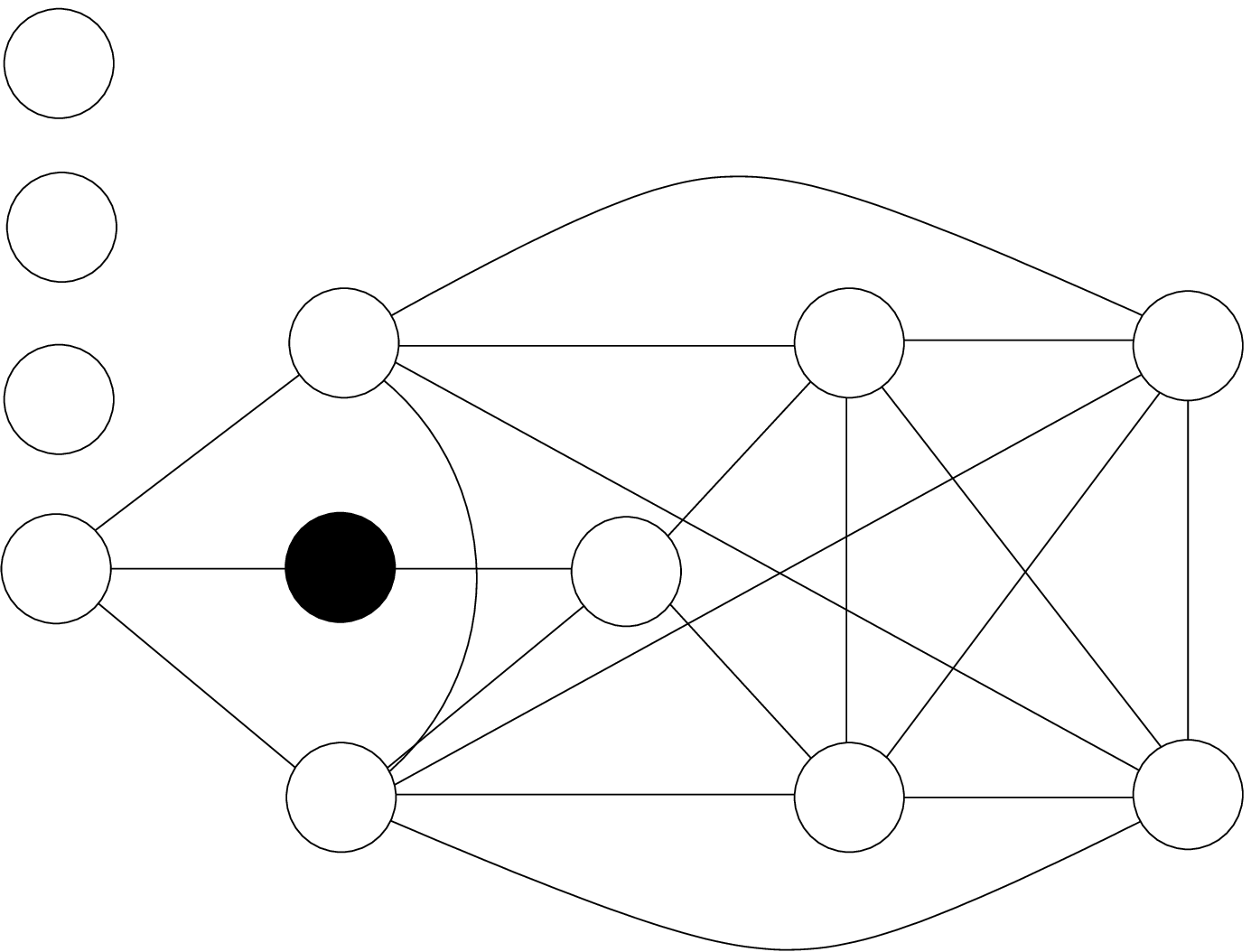, width=2.4cm} & \raisebox{10mm}[0mm][0mm]{$\rightarrow$} & \epsfig{file=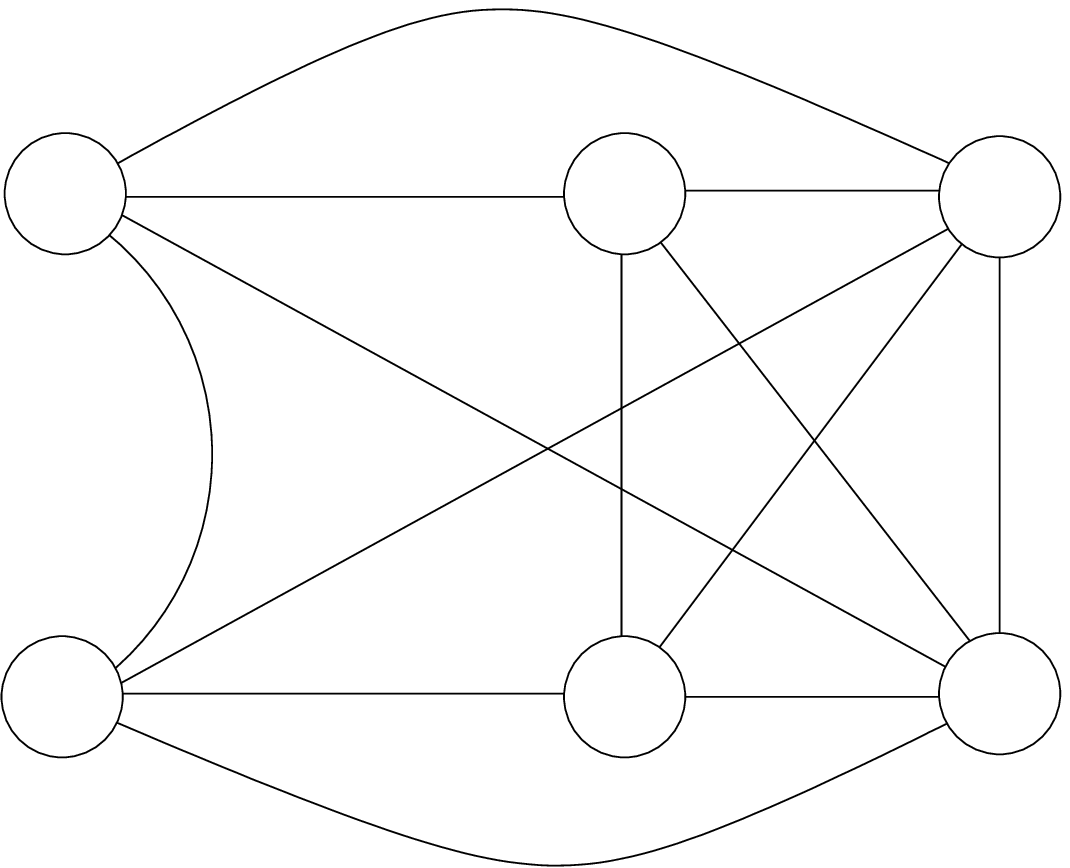, width=1.8cm}& \raisebox{10mm}[0mm][0mm]{$\xrightarrow{\T{8}\begin{smallmatrix}\mbox{choose any}\\\mbox{degree 4 vertex}\end{smallmatrix}}$}& \raisebox{9mm}[0mm][0mm]{$K_1$}\\
\mathcal G_1=G & & \mathcal G_2=3K_1+G_2 & & \mathcal G_3=G_3 & & \mathcal G_4=K_1\\
\d_1=2 & & \d_2=2, m_2=3 & & \d_3=4, m_3=0 & & \d_4=0, m_4=1
\end{array}$

Now $(m_2+1-\d_2)+(m_3+1-\d_3)=(2)+(-3)=-1$.

\item Suppose we choose $x$ as the first vertex. So we have

$\begin{array}{ccccccccc}
\epsfig{file=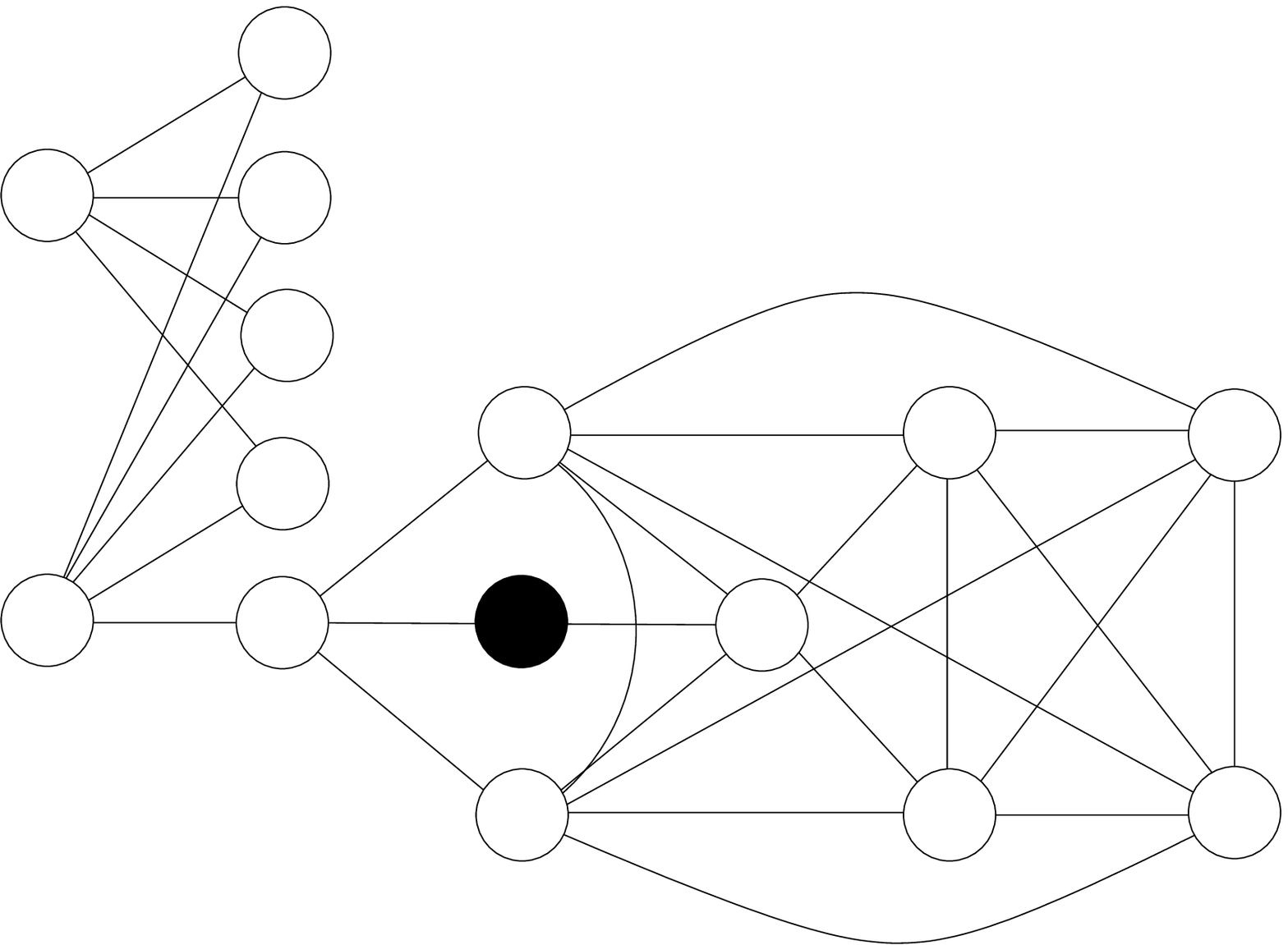, width=3cm} &\raisebox{10mm}[0mm][0mm]{$\rightarrow$}&
\epsfig{file=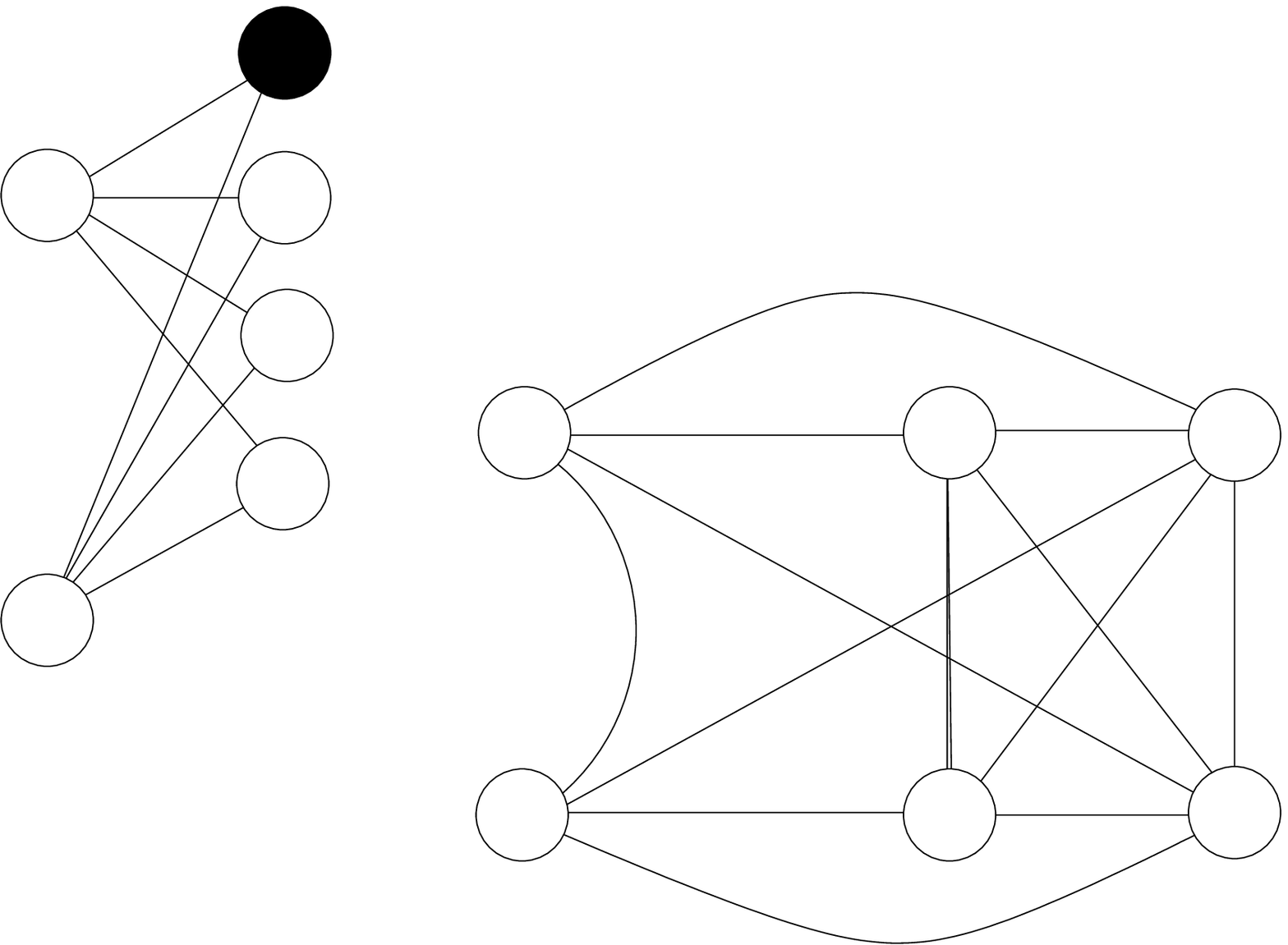, width=3cm} & \raisebox{10mm}[0mm][0mm]{$\rightarrow$} & \epsfig{file=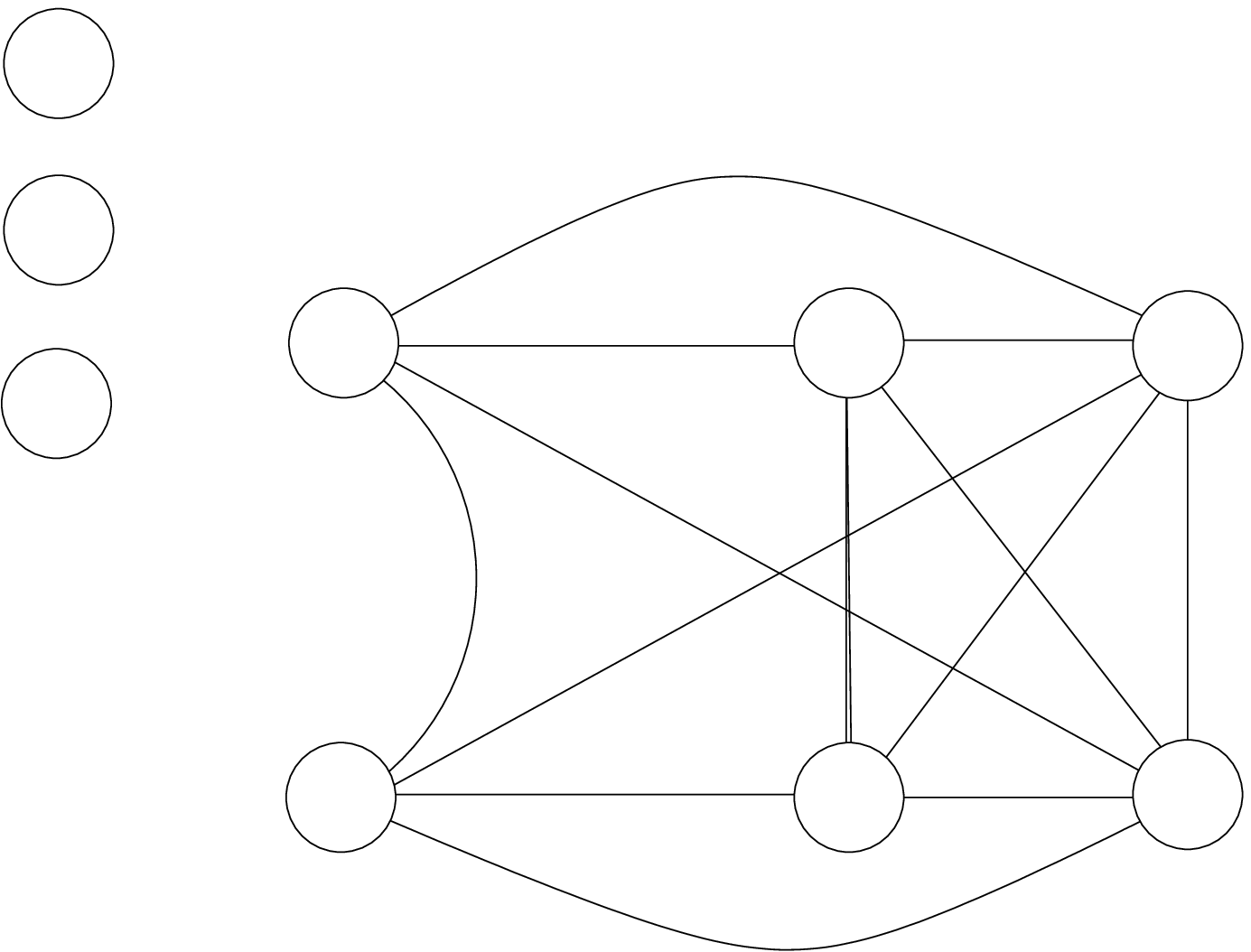, width=2.4cm}& \raisebox{10mm}[0mm][0mm]{$\xrightarrow{\T{8}\begin{smallmatrix}\mbox{choose any}\\\mbox{degree 4 vertex}\end{smallmatrix}}$} & \raisebox{9mm}[0mm][0mm]{$K_1$}\\
\mathcal G_1=G & & \mathcal G_2=G_2 & & \mathcal G_3=3K_3+G_3 & & \mathcal G_4=K_1\\
\d_1=2 & & \d_2=2, m_2=0 & & \d_3=4, m_3=3 & & \d_4=0, m_4=1
\end{array}$

Now $(m_2+1-\d_2)=-1$.
\end{enumerate}
So every $\d$-sequence of $G$ does not satisfy condition \eqref{eq-cond}.  \rsq
\end{example}

\ms
\nt The following theorem gives a sufficient condition for a graph to have minimal strength.

\begin{theorem}\label{thm-sufficient} For a graph $G$ of order $p$ with $p-2\ge \d(G)=\d_1\ge 1$, if there is a $\d$-sequence $\{\mathcal G_i\}_{i=1}^s$ of $G$ such that
\begin{equation}\label{eq-cond} \widetilde{z}_i(G)=\sum\limits_{j=2}^i \widetilde{y}_j(G)\ge 0 \mbox{ for } 2\le i\le s,\end{equation}
then $str(G)=p+\d(G)$.
\end{theorem}

\nt Note that sum with empty term is treated as zero as usual. If there is no ambiguity, we will write $\widetilde{y}_j(G)$ as $\widetilde{y}_j$ and $\widetilde{z}_i(G)$ as $\widetilde{z}_i$.

\begin{proof} Let $\{\mathcal G_i\}_{i=1}^s$ be a $\d$-sequence of $G$ satisfying condition \eqref{eq-cond}. Let $u_i$ be a vertex in $G_i$ of degree $\d_i$, which is deleted from $G_i$ to obtain $\mathcal G_{i+1}$, $1\le i\le s-1$. Now $\mathcal G_1 = G_1=G$. We shall construct a numbering $f$ of $G$ such that $str_f(G)= p+\d_1$.

\ms\nt Label $u_1$ by $p$ and all its neighbors by 1 to $\d_1$ in arbitrary order. This guarantees that the largest induced edge label is $p+\d_1$ at this stage.

\ms\nt Suppose we have labeled vertices in $V(G)\setminus V(\mathcal G_{i+1})$ by using the labels in $[1, \sum\limits_{j=1}^i\d_j]\cup [p+1-\sum\limits_{j=1}^i(m_j+1), p]$, where $1\le i\le s-1$. Moreover, the neighbors of $u_i$ are labeled by labels in $[1+\sum\limits_{j=1}^{i-1}\d_j, \sum\limits_{j=1}^i\d_j]$, and all induced edge labels are at most $p+\d_1$, up to now. Note that, sum with empty term is treated as zero.

\ms\nt Now we consider the graph $\mathcal G_{i+1}$.
\begin{enumerate}[(a)]

\item Suppose $2\le i+1<s$. We label the $m_{i+1}$ isolated vertices of $\mathcal G_{i+1}$ by labels in $[p+2-\sum\limits_{j=1}^{i+1}(m_j+1) , p-\sum\limits_{j=1}^{i}(m_j+1)]$ respectively (if $m_{i+1}=0$, then this process does not exist),  $u_{i+1}$ by $p+1-\sum\limits_{j=1}^{i+1}(m_j+1)$ and its neighbors  by label in $[1+\sum\limits_{j=1}^{i}\d_j, \sum\limits_{j=1}^{i+1}\d_j]$ respectively.

    Now, the vertices of $V(G)\setminus V(\mathcal G_{i+2})$ are labeled by using the labels in $[1, \sum\limits_{j=1}^{i+1}\d_j]\cup [p+1-\sum\limits_{j=1}^{i+1}(m_j+1), p]$.

    Since each isolated vertex of $\mathcal G_{i+1}$ is only adjacent to some neighbors of $u_i$, and $u_{i+1}$ may be adjacent with some neighbors of $u_i$, the largest new induced edge label related to these vertices is
    \begin{align*} & \quad p+1-\sum\limits_{j=1}^{i+1}(m_j+1)+\sum\limits_{j=1}^{i}\d_j
    =p+1-\sum\limits_{j=2}^{i}(m_j+1-\d_j)-(m_{i+1}+1)-(m_1+1)+\d_1\\
    & = p-z_i -m_{i+1} -1+\d_1 <p+\d_1 \tag{since $m_1=0$}.\end{align*}
    The largest new induced edge label related to $u_{i+1}$ and its neighbors in $\mathcal G_{i+1}$ is
    \begin{align*} & \quad p+1-\sum\limits_{j=1}^{i+1}(m_j+1)+\sum\limits_{j=1}^{i+1}\d_j
    =p+1-\sum\limits_{j=2}^{i+1}(m_j+1-\d_j)-(m_1+1-\d_1)\\ & =p+1-\widetilde{z}_{i+1}-(m_1+1-\d_1)\le p+1-(1-\d_1)=p+\d_1.\end{align*}
Repeat this process until $i+1=s$.
\item Suppose $i+1=s$. Now, $\mathcal G_s=m_sK_1$ with $m_s\ge 1$ or $m_sK_1+K_r$ for some $r\ge 2$ and $m_s\ge 0$. In this case, the set of unused labels is $[1+\sum\limits_{j=1}^{s-1}\d_j, p-\sum\limits_{j=1}^{s-1} (m_j+1)]$. That is, $m_s=p-\sum\limits_{j=1}^{s-1} (m_j+1+\d_j)$ or $m_s+r=p-\sum\limits_{j=1}^{s-1} (m_j+1+\d_j)$.

    When $\mathcal G_s=m_sK_1$. The process is the same as in the above case. Hence we have a numbering for $G$ with the strength $p+\d_1$.

    When $\mathcal G_s=m_sK_1+K_r$, where $\d_s+1=r=p-m_s-\sum\limits_{j=1}^{s-1} (m_j+1+\d_j)$. We label the $m_{s}$ isolated vertices of $\mathcal G_{s}$ by labels in $[p-m_s+1-\sum\limits_{j=1}^{s-1}(m_j+1) , p-\sum\limits_{j=1}^{s-1}(m_j+1)]$ respectively (if $m_{s}=0$, then this process does not perform). Finally, label the vertices of $K_r$ by labels in $[1+\sum\limits_{j=1}^{s-1}\d_j,p-m_s-\sum\limits_{j=1}^{s-1}(m_j+1)]$, respectively. Then the largest new induced edge labels related to the neighbors of $u_{s-1}$ is
    \[p-\sum_{j=1}^{s-1}(m_j+1)+\sum_{j=1}^{s-1}\d_j=p-(m_1+1-\d_1)-\widetilde{z}_{s-1}\le p-1+\d_1<p+\d_1.\]

 The largest new induced edge labels in $K_r$ is
\begin{align*} & \quad [p-m_s-\sum\limits_{j=1}^{s-1} (m_j+1)]+[p-1-m_s-\sum\limits_{j=1}^{s-1} (m_j+1)]\\ &=[p-m_s-\sum\limits_{j=1}^{s-1} (m_j+1)]+[\sum\limits_{j=1}^s \d_j]=p+1-\sum\limits_{j=1}^{s} (m_j+1-\d_j)\\ &\le p+1-m_1-1+\d_1=p+\d_1.\end{align*}
\end{enumerate}
Hence we have a numbering $f$ such that $str_f(G)=p+\d_1$. Therefore, $str(G)\le p+\d_1$. By Lemma~\ref{lem-d>0}, $str(G)=p+\d_1$.
\end{proof}

\begin{example}\label{ex-counterexample1} Consider the graph $G_1$ described in Example~\ref{ex-example}. Using the second $\d$-sequence of $G_1$ and following the construction in the proof of Theorem~\ref{thm-sufficient} we have the following strength labeling $f$ of $G_1$ such that $str_f(G_1)=14$.

\centerline{\epsfig{file=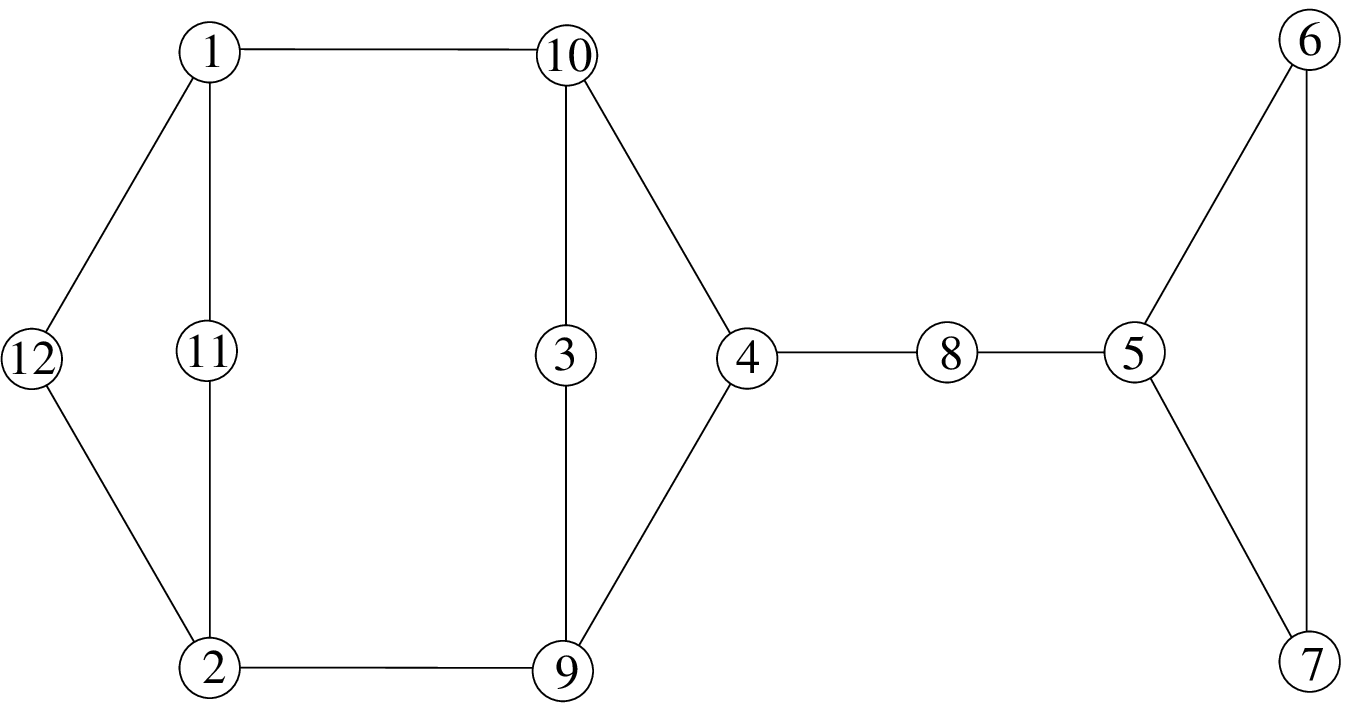, width=7cm}} \rsq
\end{example}

\nt From Example~\ref{ex-counterexample}, $G$ does not satisfy the hypothesis of Theorem~\ref{thm-sufficient}, but there is a strength labeling $f$ for it with $str_f(G)=p+\d(G)=17$ as follows:

\centerline{\epsfig{file=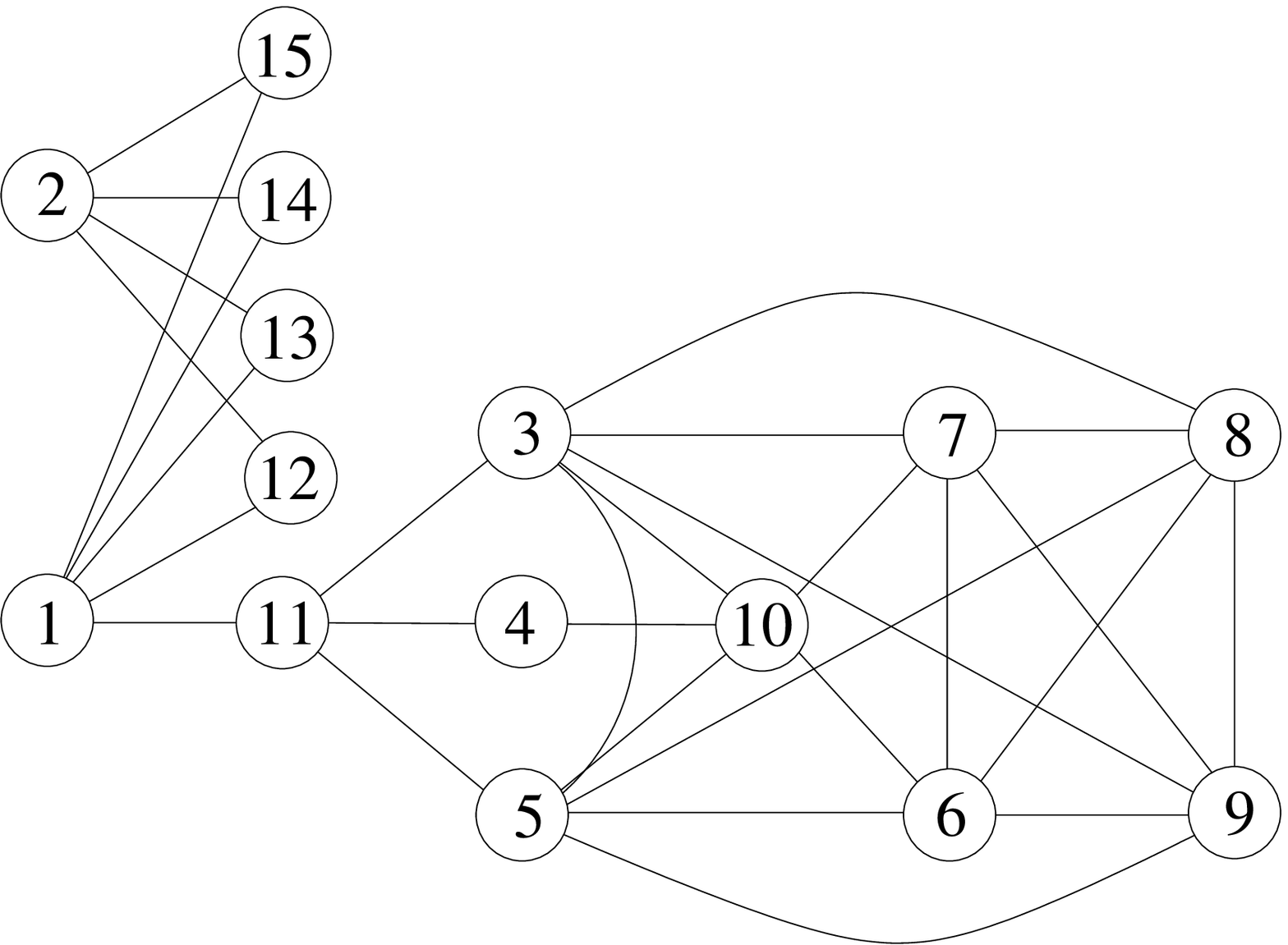, width=5cm}}

\nt So the converse of Theorem~\ref{thm-sufficient} is not true.

\ms

\nt Thus, Theorem~\ref{thm-sufficient} provides a solution to Problem~\ref{pbm-minstr}. Note that every tree $T$ has the property that $\d_i=1$ and $m_i\ge 0$ for each $i\ge 2$. We immediately have $str(T)=|V(T)|+1$ and the following corollary that answers more than what Problem~\ref{pbm-lobster} asks.

\begin{corollary} If $G$ is a forest without isolated vertex, then $str(G) = |V(G)| + 1$. \end{corollary} 

\begin{corollary} The one-point union of cycles $G$ of order $p$ has $str(G)=p+2$. \end{corollary}

\begin{proof} Remove a degree 2 vertex that is adjacent to the maximum degree vertex of $G$ and its neighbors to obtain a subgraph $\mathcal G_2$, which is a disjoint union of path(s). So, $G$ admits a $\d$-sequence that satisfies~\eqref{eq-cond}.  \end{proof}

\begin{corollary} If $G$ is a wheel or fan graph of order $p$, then $str(G)=p+\d(G)$. \end{corollary}

\ms

\nt In constructing a $\d$-sequence of $G$, if we change the choice of choosing a vertex of degree $\d_i$ to a vertex of degree $d_i$, then we get another sequence of subgraphs of $G$. This sequence is called a {\it $d$-sequence} of $G$. Let $y_j(G)=m_j+1-d_j$, $z_i(G)=\sum\limits_{j=2}^i y_j(G)$ and denote $d_1$ by $d_G$. By the same argument as in proving Theorem~\ref{thm-sufficient} we have
\begin{theorem}\label{thm-sufficient2} For a graph $G$ of order $p$ with $p-2\ge \d(G)\ge 1$, if there is a $d$-sequence $\{\mathcal G_i\}_{i=1}^s$ of $G$ such that
\begin{equation}\label{eq-cond2} z_i(G)\ge 0 \mbox{ for } 2\le i\le s,\end{equation}
then $str(G)\le p+d_G$. The equality holds if $d_G = \d(G)$.
\end{theorem}



\begin{example} Consider the graph $G$ in Example~\ref{ex-counterexample}. Following is a $d$-sequence of $G$.

$\begin{array}{ccccccc}
\epsfig{file=strength-1a.eps, width=3cm} &\raisebox{10mm}[0mm][0mm]{$\rightarrow$}&
\epsfig{file=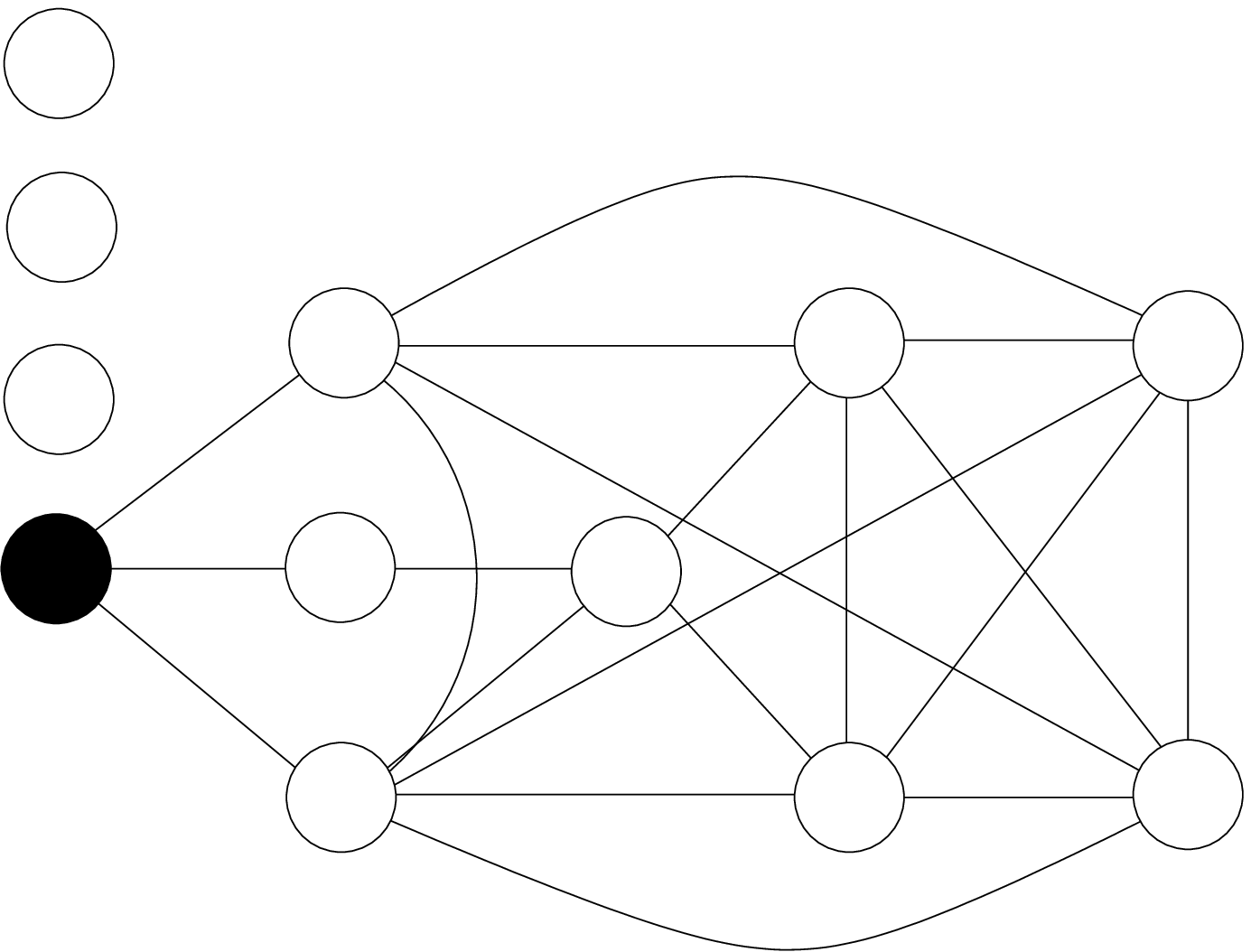, width=2.4cm} & \raisebox{10mm}[0mm][0mm]{$\rightarrow$} & \raisebox{5mm}[0mm][0mm]{\epsfig{file=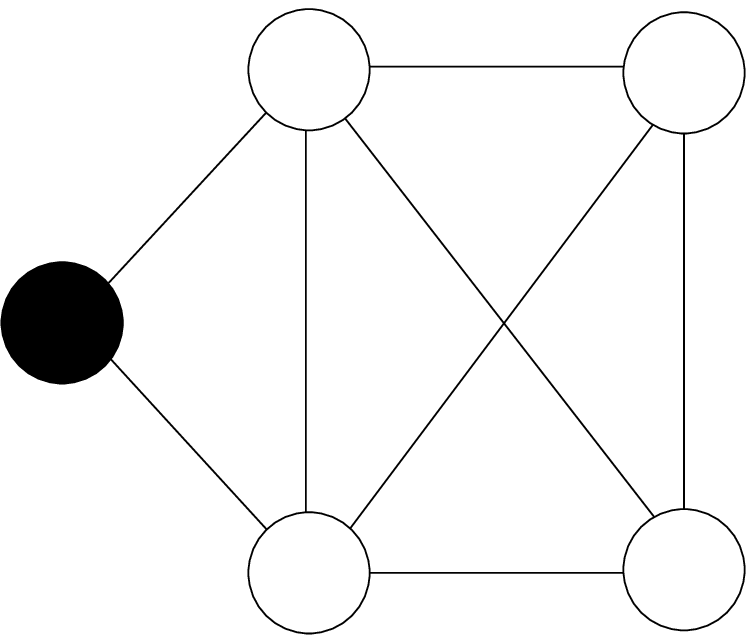, width=1.3cm}}&
\raisebox{10mm}[0mm][0mm]{$\rightarrow$} & \raisebox{5mm}[0mm][0mm]{\epsfig{file=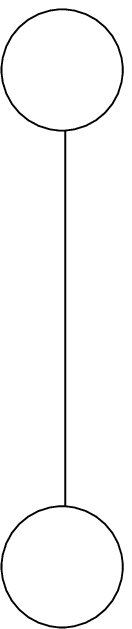, width=2mm}}\\
\mathcal G_1=G & & \mathcal G_2=3K_1+G_2 & & \mathcal G_3=G_3 & & \mathcal G_4=K_2\\
d_G=d_1=\d(G)=2 & & d_2=3, m_2=3 & & d_3=2, m_3=0 & & d_4=1, m_4=0
\end{array}$

Now $z_2=1$, $z_3=0$ and $z_4=0$. By Theorem~\ref{thm-sufficient2}, $str(G)=17$ as shown in Example~\ref{ex-counterexample1}. \rsq
\end{example}

\begin{lemma}\label{lem-forest} Let $T$ be a forest without isolated vertex of order at least $3$ and let $P_T$ be the set of pendant vertices that are adjacent to a vertex of degree at least $2$. There is a $\d$-sequence of $T$ of length $s$ such that $\widetilde{z}_s(T)\ge |P_T|-|N_T(P_T)|$, where $\widetilde{z}_i(T)$ is defined in \eqref{eq-cond}. Moreover, all $\widetilde{z}_i(T)$ satisfy \eqref{eq-cond}.
\end{lemma}
\begin{proof}
\nt Obviously, the lemma holds when $T$ is of order $3$. Suppose the lemma holds when the order of $T$ is $k$ or less, where $k\ge 3$.

\ms\nt Now consider a forest $T$ of order $k+1$. Choose a vertex $u\in P_T$. Let $v$ be the vertex adjacent to $u$ with degree $d$. We shall consider the forest $T-u-v$.

\ms\nt Suppose $T=K_{1,k}$ which is a star, then $T-u-v=(k-1)K_1$ and $\{T, (k-1)K_1\}$ is a $\d$-sequence of $T$. Note that $|P_T|=k$ and $|N_T(P_T)|=1$. Clearly $\widetilde{z}_2(T)=(k-1)+1-0=k>|P_T|-|N_T(P_T)|$.

\ms\nt Now we assume that $T$ is not a star. Let $T-u-v = mK_1+T'$, where $m\ge 0$ and $T'$ is a forest without isolated vertex.

\ms\nt Suppose the order of $T'$ is 2, then $T'=K_2$ and $\{T, mK_1+K_2\}$ is a $\d$-sequence of $T$, where $m=k-3$. Now $\widetilde{z}_2(T)=m=k-3$. If $T=K_{1,k-2}+K_2$, then $|P_T|=k-2$ and $|N_T(P_T)|=1$.  We get  $\widetilde{z}_2(T)=|P_T|-|N_T(P_T)|$. If $T$ is a tree, then $|P_T|=k-1$ and $|N_T(P_T)|=2$.  So we still get $\widetilde{z}_2(T)=|P_T|-|N_T(P_T)|$.

\ms\nt Suppose the order of $T'$ is greater than 2. By induction assumption there is a $\d$-sequence of $T'$, say $\{T'=\mathcal T_1,\mathcal T_2,\dots, \mathcal T_s\}$, such that $\widetilde{z}_s(T')\ge |P_{T'}|-|N_{T'}(P_{T'})|$ for some $s\ge 2$.

\ms\nt Now, consider $\{T, \mathcal T_1,\dots, \mathcal T_s\}$ of $T$. Note that $\widetilde{z}_2(T) =m+1-\d(T')=m$. Let $a$ be the number of vertices of degree 2 in $T$ but degree 1 in $T'$, then $|P_{T'}|=|P_T|-(m+1)+a$, where $a\ge 0$. Let $b=|N_{T'}(P_{T'})\setminus N_T(P_T)|$.
Therefore, $|N_{T'}(P_{T'})|\le |N_T(P_T)|-1+b$ (some vertices in $N_T(P_T)$ may not be in $N_{T'}(P_{T'})$). 

\ms\nt If a vertex $w$ in $T$ of degree 2 becomes of degree 1 in $T'$, then $w$ may be in $P_{T'}$ and may have at most one neighbor in $N_{T'}(P_{T'})$ which is not in $N_T(P_T)$. So $a\ge b$.

\nt Now, we have 
\begin{align*}
\widetilde{z}_{s+1}(T)& =\widetilde{z}_2(T)+\widetilde{z}_s(T')=m+\widetilde{z}_s(T')\ge
m+|P_{T'}|-|N_{T'}(P_{T'})|\\
& \ge m+[|P_T|-(m+1)+a]-[|N_T(P_T)|-1+b]\\ &= |P_T|-|N_T(P_T)|+a-b\ge |P_T|-|N_T(P_T)|.
\end{align*}
By induction the lemma holds for any forest $T$ of order greater than 2.
\end{proof}

\begin{rem}\label{rem-star}
In the proof of Lemma~\ref{lem-forest}, we can see that $\widetilde{z}_2(T)=|P_T|-|N_T(P_T)|+1$, if $T$ is a star.\end{rem} 

\begin{theorem}\label{thm-extension1} Keep all notation defined in Theorem~\ref{thm-sufficient} and Lemma~\ref{lem-forest}. Let $H$ be a graph with $\d(H)\ge 1$. Suppose $\{H=\mathcal H_1, \dots, \mathcal H_s\}$ is a $d$-sequence of $H$. Suppose
\begin{equation}\label{eq-Z1} Z=\min\{z_i(H)\; |\; 2\le i\le s\}\end{equation}
is not positive.
Suppose $T$ is a graph with a $d$-sequence $\{T=\mathcal T_1, \dots, \mathcal T_t\}$  satisfying \eqref{eq-cond2}.
Let $G = H + T$. If $z_t(T)\ge d_H-Z$, then $|V(G)| + \d(G)\le str(G) \le |V(G)| + d_T$. \end{theorem}

\begin{proof} Since $\mathcal T_t$ is either $mK_1$ for some $m\ge 1$ or $mK_1+K_r$ for some $m\ge 0$, we have the following two cases.

\nt For the first case, $y_t(T)=m+1$ and $\{\mathcal T_1, \dots, \mathcal T_{t-1}, mK_1+H, \mathcal H_2, \dots, \mathcal H_s\}$ is a $d$-sequence of $G$. Then\\
$z_j(G)=z_j(T)\ge 0$, $2\le j\le t-1$;\\
$z_t(G)=z_{t-1}(T)+y_t(H+T)=z_{t-1}(T)+[m+1-d_H]=z_t(T)-d_H\ge z_t(T)-d_H+Z\ge 0$;\\ $z_{t+1}(G)=z_{t-1}(T)+(m+1-d_H)+z_2(H)=z_t(T)-d_H+y_2(H)=z_t(T)-d_H+z_2(H)$.\\
In general, \begin{align*}z_{t+j}(G)& =z_{t-1}(T)+(m+1-d_H)+z_{j+1}(H)=z_t(T)-d_H+z_{j+1}(H)\\ & \ge z_t(T)-d_H+Z\ge 0, \quad 1\le j\le s-1.\end{align*}

\nt For the last case, $\{\mathcal T_1, \dots, \mathcal T_t, \mathcal H_1, \mathcal H_2, \dots, \mathcal H_s\}$ is a $d$-sequence of $G$. We will get that\\
$z_j(G)=z_j(T)\ge 0$, $2\le j\le t$;\\
$z_{t+1}(G)=z_t(G)+[1-d_H]=z_t(T)+1-d_H>0$;\\
$z_{t+j}(G) =z_t(T)+1-d_H+z_{j}(H)>0$, $2\le j\le s$.

\ms\nt By Theorem~\ref{thm-sufficient2}, $str(G)\le |V(G)|+d_T$. The lower bound follows from Lemma~\ref{lem-d>0}.
\end{proof}

\begin{rem}\label{rem-connect} In Theorem~\ref{thm-extension1}, suppose $H$ and $T$ are connected. Let $x\in V(\mathcal T_{t-1})$ be chosen to construct $\mathcal T_t$, $v\in V(\mathcal T_{t-1})$ be a neighbor of $x$ and let $u\in V(H)$ which is not chosen to construct $\mathcal H_2$. We add an edge $vu$ to the graph $H+T$. All the $y_j$ values of this connected graph are the same as those of $H+T$.
\end{rem}

\begin{theorem} Keep all notation defined in Theorem~\ref{thm-sufficient} and Lemma~\ref{lem-forest}. Let $H$ be a graph with $\d(H)\ge 1$. Suppose $\{H=\mathcal H_1, \dots, \mathcal H_s\}$ is a $d$-sequence of $H$. Suppose
\[Z=\min\{z_i(H)\; |\; 2\le i\le s\}\le 0.\]
Let $G = H + T$, where $T$ is a forest without isolated vertex of order at least $3$. If $|P_T|-|N_T(P_T)|\ge d_H-Z$, then $str(G) = |V(G)| + 1$. \end{theorem}

\begin{proof} By Lemma~\ref{lem-forest} there is a $\d$-sequence $\{\mathcal T_1, \dots, \mathcal T_t\}$ of $T$ such that $\widetilde{z}_t(T)\ge |P_T|-|N_T(P_T)|$ and all $\widetilde{z}_i(T)$ satisfy \eqref{eq-cond}. Since this $\d$-sequence is a particular $d$-sequence of $T$, it satisfies the condition of Theorem~\ref{thm-extension1}. Since $d_T=1$ now, by Theorem~\ref{thm-extension1} we have $str(G)\le |V(G)|+1$.
Hence we have the theorem since $\d(G)=1$.
\end{proof}

\begin{corollary} Let $H$ be a graph with $\d(H)\ge 1$. Suppose $\{H=\mathcal H_1, \dots, \mathcal H_s\}$ is a $d$-sequence of $H$. Suppose
\[Z=\min\{z_i(H)\; |\; 2\le i\le s\}\le 0.\]
Let $G = H + K_{1,k}$. If $k\ge d_H-Z$, then $str(G) = |V(G)| + 1$. \end{corollary}

\nt Let $G=K_{m,n}$ for $n\ge m\ge 1$. It is proven in\cite[Theorem 3.5]{Ichishima+MB+Oshima18} that $str(G) = |V(G)| + m$. So we have

\begin{corollary} There exists a graph $G$ with $str(G)=|V(G)|+\d(G)$ for each $\d(G)\ge 1$. \end{corollary}


\nt Suppose $H$ is a graph with $str(H) > |V(H)| + \d(H)$. So there is a $d$-sequence of $H$ with $d_H=\d(H)$ such that $Z$ is non-positive, where $Z$ is defined in \eqref{eq-Z1}. Let $T=K_{m,n}$. Since $\{T, (n-1)K_1\}$ is a $d$-sequence of $T$, $z_2(T)=n$. Suppose $n\ge d_H - Z$ and $m=\d(H)$. Since $n\ge d_H - Z\ge \d(H)$, $\d(H+T)=\d(H)$. By Theorem~\ref{thm-extension1} we have $str(H+T)=|V(H+T)|+\d(H)$. So, together with Remark~\ref{rem-connect}, we have the following theorem.

\begin{theorem}\label{thm-extension2} For every graph $H$, either $str(H) = |V(H)| + \d(H)$ or $H$ is a proper subgraph of a graph $G$ such that $str(G) = |V(G)| + \d(G)$ with $\d(G) = \d(H)$. \end{theorem}


\begin{example} It is easy to obtain
a $\d$-sequence of $Q_4$ with
\begin{align*}
\d_1 &= 4;\\
\d_2& = 2, m_2 = 0,   \mbox{ so } \widetilde{z}_2 = -1;\\
\d_3 &= 1, m_3 = 0,  \mbox{ so }  \widetilde{z}_3 = -1;\\
\d_4 &= 1, m_4 = 1,   \mbox{ so }  \widetilde{z}_4 = 0;\\
\d_5& = 0, m_5 = 3,  \mbox{ so } \widetilde{z}_5 = 4.\end{align*}
Thus, $Z = -1$. Let $n= 4 + 1 = 5$. By the construction before Theorem~\ref{thm-extension2}, we have  $str(Q_4+K_{4,5}) = 25 + 4 = 29$. A required labeling can be obtained by similarly following the proof of Theorem~\ref{thm-sufficient}. Moreover, adding an edge joining a vertex of degree 5 of $K_{4,5}$ and a vertex of $Q_4$ gives a connected graph $G$ that contains $Q_4$ as a proper subgraph with $str(G) = str(Q_4+K_{4,5})$ as required. \rsq
\end{example}

\begin{example}For $2$-regular graphs $\mathcal{C}_k$ with exactly $k\ge 1$ odd cycles, we have $Z = -k+1$. Let $n =k+1$. By the construction before Theorem~\ref{thm-extension2}, we have $str(\mathcal{C}_k+K_{2,k+1}) = |V(\mathcal{C}_k)|+ (k+3)+ 2$. \rsq
\end{example}

\section{New Lower Bounds}

\begin{theorem}\label{thm-indep-no} Suppose $G$ is a graph of order $p$ with independent number $\alpha$, then $str(G)\ge 2p-2\alpha+1$.
\end{theorem}
\begin{proof} For any numbering for $G$, by pigeonhole principle, at least two labels in $[p-\alpha, p]$ are assigned to two adjacent vertices. So the induced edge label at least $2p-2\alpha+1$. This completes the proof.
\end{proof}

\begin{corollary}\label{cor-necessary-indep-no} Suppose $G$ is a graph of order $p$ with minimum degree $\d$. Suppose $str(G)=p+\d$, then $\alpha \ge \left\lceil\frac{p-\d+1}{2}\right\rceil$, where $\alpha$ is the independence number of $G$.
\end{corollary}
\begin{proof} From Theorem~\ref{thm-indep-no}, we have $\alpha\ge \frac{p-\d+1}{2}$.
\end{proof}

\nt Let $G$ be a graph of order $p$. Let \begin{align*}x_i&=\min\left\{| N_G(S)\setminus S|\;:\; |S|=i\right\};\\
\xi=\xi(G) &=\max\{x_i-i+1\;|\; 1\le i\le p-1\}.\end{align*}

\begin{theorem}\label{thm-neighborhood}
Let $G$ be a graph of order $p$, then $str(G)\ge p+\xi$.
\end{theorem}
\begin{proof} Let $\xi=x_i-i+1$ for some $i$. Let $f$ be a strength labeling of $G$. Consider the labels in $[p-i+1, p]$. Let $T=f^{-1}([p-i+1, p])$, then $|T|=i$. Now $|f(N_G(T)\setminus T)|=|N_G(T)\setminus T|\ge x_i$. Let $a$ be the largest label in $f(N_G(T)\setminus T)$. There is a vertex $u\in N_G(T)\setminus T$ such that $f(u)=a$. Moreover, $u$ is adjacent to $v\in T$. Thus,
\[str_f(G)\ge f(v)+f(u)\ge p-i+1+a\ge p-i+1+x_i=p+\xi.\]
\end{proof}

\nt Thus, we provided 2 good bounds for the strength of a graph as raised in Problem~\ref{pbm-gb}. Note that Lemma~\ref{lem-d>0} is a corollary of  Theorem~\ref{thm-neighborhood} when $\xi = \d = x_1 \ge x_i-i+1$ for $i\ge 2$.






\begin{theorem}\label{thm-2regular} If $G=\sum_{i=1}^h C_{2m_i}+\sum_{j=1}^k C_{2n_j+1}$, where $m_i\ge 2$, $n_j\ge 1$ and $h+k\ge 1$, then $str(G)=\max\{p+2, p+1+k\}$. \end{theorem}

\begin{proof} Note that if $h=0$, then the first summand does not appear, similarly for the second summand. Now, $\alpha(G)=\sum_{i=1}^{h} m_i+\sum_{j=1}^{k} n_j$. By Theorem~\ref{thm-indep-no}, we have $str(G)\ge p+1+k$.

\ms
\nt Let $H=\sum_{i=1}^h C_{2m_i}$ and $K=\sum_{j=1}^k C_{2n_j+1}$. Let $M=\sum_{i=1}^h m_i$ such that $M=0$ when $h=0$.

\ms
\nt We shall construct a numbering $f$ on $G$. If $h\ge 1$, we first label $H$ by $[1, M]\cup [p-M+1, p]$ as follows:

\ms\nt Label $C_{2m_1}$ by $1, p, 2, \dots, m_1,  p-m_1+1$ in the natural order. In general, for $i\ge 2$, we label the vertices of even cycle $C_{2m_i}$ by
$1+\sum_{l=1}^{i-1} m_l$, $p-\sum_{l=1}^{j-1}m_l$, $2+\sum_{l=1}^{j-1} m_l$, $\dots$, $\sum_{l=1}^{i}m_l$, $p-\sum_{l=1}^{j}m_l+1$ in the natural order. Continue this process until $i=h$. Hence the maximum induced edge label is $p+2$.

\ms\nt If $k\ge 1$, then we label the vertices of odd cycle $C_{2n_1+1}$ by $M+1, p-M, M+2, \dots, p-M-n_1+1, M+n_1+1$ in the natural order. Up to now, the maximum induced edge label is still $p+2$.

\ms\nt Now we label the vertices of odd cycle $C_{2n_2+1}$ by $M+n_1+2$, $p-M-n_1$, $M+n_1+3$, $\dots$, $p-M-n_1-n_2+1$, $M+n_1+n_2+2$ in the natural order. Note that $(M+n_1+2)+(M+n_1+n_2+2)=2M+(2n_1+1)+n_2+2\le p$. So, the current maximum induced edge label is $p+3$.

\ms\nt In general, for $j\ge 2$, we label the vertices of odd cycle $C_{2n_j+1}$ by $M+1+\sum_{l=1}^{j-1} (n_l+1)$, $p-M-\sum_{l=1}^{j-1}n_l$, $M+2+\sum_{l=1}^{j-1} (n_l+1)$, $\dots$, $p-M-\sum_{l=1}^{j}n_l+1$, $M+\sum_{l=1}^{j}(n_l+1)$ in the natural order. Note that $(M+1+\sum_{l=1}^{j-1} (n_l+1))+(M+\sum_{l=1}^{j} (n_l+1))=2M+\sum_{l=1}^{j-1}(2n_l+1)+n_j+j\le p-2+j$. So, the current maximum induced edge label is $p+j+1$.

\ms\nt Continue this process until $j=k$. Hence we have $str_f(G)=p+k+1$.
\end{proof}

\begin{example} Consider $G=C_4+C_6+C_5+C_5+C_7$. Now, $p=27$ and $k=3$.

\nt We label\\
$C_4$ by $[1, 27, 2, 26]$ (max. induced edge label is 29);\\
$C_6$ by $[3, 25, 4, 24, 5, 23]$ (max. induced edge label is 29);\\
$C_5$ by $[6, 22, 7, 21, 8]$; (max. induced edge label is 29);\\
$C_5$ by $[9, 20, 10, 19, 11]$ (max. induced edge label is 30);\\
$C_7$ by $[12, 18, 13, 17, 14, 16, 15]$ (max. induced edge label is 31).

\nt So $str(G)=31$. \rsq
\end{example}



\nt Let $G\times H$ be the Cartesian product of graphs $G$ and $H$.

\begin{lemma}\label{lem-bipartite} Let $G$ be a bipartite graph with bipartition $(X,Y)$ such that $|X|=|Y|=m$. Suppose there is a numbering $f$ of $G$ such that $f(X)=[1,m]$, then there is a numbering $F$ of $G\times K_2$ such that $F([1,2m])=\widetilde{X}$, where $(\widetilde{X},\widetilde{Y})$ is a bipartition of $G\times K_2$. Moreover, $str_F(G\times K_2)=5m+1$.
\end{lemma}

\nt Note that the following proof is modified from the proof of Theorem~3.10 in~\cite{Ichishima+MB+Oshima18}.
\begin{proof} Note that, from the hypotheses, $f(Y)=[m+1, 2m]$. Let $u$ and $v$ be vertices of $K_2$. Then\\ $\widetilde{X}=\{(x, u)\;|\; x\in X\}\cup \{(y, v)\;|\; y\in Y\}$ and $\widetilde{Y}=\{(y, u)\;|\; y\in Y\}\cup \{(x, v)\;|\; x\in X\}$.

\nt Define $F: V(G\times K_2)\to [1,4m]$ by
\begin{align*}
F(x,u) & = f(x)\in[1,m]\\
F(y,v) & = (3m+1)-f(y)\in[m+1,2m]\\
F(x,v) & = (3m+1)-f(x)\in[2m+1,3m]\\
F(y,u) & = 2m+f(y)\in[3m+1,4m]
\end{align*}
Clearly $F(\widetilde{X})=[1, 2m]$.

\ms\nt Now $F(x, u)+F(x,v)=f(x)+(3m+1)-f(x)=3m+1$ and $F(y, u)+F(y,v)=2m+f(y)+(3m+1)-f(y)=5m+1$. Suppose $(x_1, u)$ and $(x_2,u)$ are adjacent in $G\times K_2$. By definition, $F(x_1, u)+F(x_2,u)\le 2m$.
Similarly if $(y_1, v)$ and $(y_2,v)$ are adjacent in $G\times K_2$, then $F(y_1, v)+F(y_2,v)\le 4m$. Thus $str_F(G\times K_2)=5m+1$.
\end{proof}

\nt Let $Q_n$ be the hypercube of dimension $n$, $n\ge 2$. Since there is a strength numbering $f$ of $Q_2$ satisfying the hypotheses of Lemma~\ref{lem-bipartite}, applying this lemma repeatedly we get that

\begin{theorem}\label{thm-Qn} For $n\ge 2$, $str(Q_n)\le 2^n+2^{n-2}+1$. \end{theorem}

\nt This is a known result in \cite[Theorem~3.10]{Ichishima+MB+Oshima18}.

\ms\nt
We shall improve the lower bound of the strength of $Q_n$. The vertices of $Q_n$ often used the elements of the vector space $\Z_2^n$ over $\Z_2$. Two vertices $u$ and $v$ are adjacent if and only if $u+v=e_i$, where $e_i$ is the standard basis of $\Z_2^n$. Note that, $v= -v$ for any vector $v\in \Z_2^n$. In the proofs of the following lemmas, all algebra involving vectors are over $\Z_2$.

\ms\nt
For any vertex $v$, we let $N_G[v]=N_G(v)\cup \{v\}$, the closed neighborhood of $v$.
Hence, for any subset of vertices $S$, $N_G(S)\setminus S=\left(\bigcup\limits_{v\in S}N_G[v]\right)\setminus S$. We shall omit the subscript $G$ when there is no ambiguity.

\begin{lemma}\label{lem-2neighbor}
If $u$ and $v$ be two distinct vertices of $Q_n$, $n\ge 3$, then $|N[u]\cap N[v]|$ is either $0$ or $2$.
\end{lemma}

\begin{proof} Suppose $u$ and $v$ are adjacent. Clearly  $|N[u]\cap N[v]| = 2$. Suppose $u$ and $v$ are not adjacent. If $z\in N[u]\cap N[v]$, then the distance between $u$ and $v$ is two. Hence $u\notin N(v)$, $v\notin N(u)$ and $u+v=e_i+e_j$ (equivalently $u+e_i=v+e_j$), where $i\ne j$. Since $z\in N(u)\cap N(v)$, $z=u+e_k=v+ e_l$ for some $k,l$. So $u+e_k+v+ e_l=\mb 0$ or $u+v=e_k+e_l$. Thus $\{i,j\}=\{k,l\}$. Hence $|N[u]\cap N[v]|=|\{u+e_i, u+e_j\}|=2$.
\end{proof}

\begin{lemma}\label{lem-3neighbor}
For any three distinct vertices $u$, $v$ and $w$ of $Q_n$, $n\ge 3$, \[|N[u]\cap N[v]\cap N[w]|=\begin{cases}0 & \mbox{if at least one of $|N[u]\cap N[v]|$, $|N[u]\cap N[w]|$ and $|N[v]\cap N[w]|$ is $0$};\\ 1 & \mbox{if all of $|N[u]\cap N[v]|$, $|N[u]\cap N[w]|$ and $|N[v]\cap N[w]|$ are $2$.}\end{cases}\]
\end{lemma}

\begin{proof} If one of $|N[u]\cap N[v]|$, $|N[u]\cap N[w]|$ and $|N[v]\cap N[w]|$ is $0$, then $|N[u]\cap N[v]\cap N[w]|=0$. Otherwise, Lemma~\ref{lem-2neighbor} implies that $|N[u]\cap N[v]|=|N[u]\cap N[w]|=|N[v]\cap N[w]|=2$.
\begin{enumerate}[(1).]
\item Suppose only one pair of $u$, $v$ and $w$ are adjacent, say $uv$ is an edge, then the distances from $w$ to $u$, and to $v$ are 2. This creates a 5-cycle which is impossible.
\item Suppose two pairs of $u$, $v$ and $w$ are adjacent, say $uv$ and $uw$ are edges. Note that $v$ and $w$ cannot be adjacent. Then $N[u]\cap N[v]=\{u,v\}$ and $N[u]\cap N[w]=\{u,w\}$. Hence $u\in N[v]\cap N[w]$. This implies that $N[u]\cap N[v]\cap N[w]=\{u\}$.

\item Suppose none of $u$, $v$ and $w$ are adjacent.  By the proof of Lemma~\ref{lem-2neighbor} we have
$u+v=e_{i_1}+e_{j_1}$ and $v+w=e_{i_2}+e_{j_2}$ for some $i_1, i_2, j_1, j_2$, $i_1\ne j_1$ and $i_2\ne j_2$. This implies that $u+w=e_{i_1}+e_{j_1}+e_{i_2}+e_{j_2}$. Since the distance of $u$ and $w$ is 2, $|\{i_1, j_1\}\cap \{i_2, j_2\}|=1$. Without loss of generality, we may assume that $i_1=i_2$. Now, $N[u]\cap N[v]=\{u+e_{i_1}, u+e_{j_1}\}$, $N[u]\cap N[w]=\{u+e_{j_1}, u+e_{j_2}\}$ and $N[v]\cap N[w]=\{v+e_{i_2}, v+e_{j_2}\}$. Here $v+e_{i_2}=v+e_{i_1}=u+e_{j_1}$.
Hence $u+e_{j_1}\in N[u]\cap N[v]\cap N[w]$. Since $u+e_{j_2}\notin N[u]\cap N[v]$, $N[u]\cap N[v]\cap N[w]=\{u+e_{j_1}\}$.
\end{enumerate}
This completes the proof. \end{proof}

\begin{theorem}\label{thm-Qn-LB} For the hypercube $Q_n$, $n\ge 2$, we have
$str(Q_2)\ge 6$;
$str(Q_3)\ge 11$; $str(Q_4)\ge 21$;
and $str(Q_n)\ge 2^n+4n-12$ for $n\ge 5$.
\end{theorem}

\begin{proof}
Keeping the notation defined in Theorem~\ref{thm-neighborhood}, we want to compute $x_i$ and $\xi$. Clearly $x_1=\d=n$.

\nt Suppose $S=\{u,v\}$ with $u\ne v$.
\begin{align*}
|N[u]\cup N[v]|=|N[u]|+ |N[v]|-|N[u]\cap N[v]|\ge 2(n+1)-2=2n.
\end{align*}

\nt So $|N(S)\setminus S|=|N[u]\cup N[v]|-|S|\ge 2n-2$. Actually when $S=\{\mb 0, e_1+e_2\}$, $|N(S)\setminus S|=2n-2$. Thus $x_2=2n-2$.

\ms

\nt Suppose $S=\{u,v, w\}$, where $u,v,w$ are distinct.
\begin{enumerate}[(1).]
\item If all $|N[u]\cap N[v]|$, $|N[u]\cap N[w]|$ and $|N[v]\cap N[w]|$ are not zero, then by Lemma~\ref{lem-3neighbor},  $|N[u]\cap N[v]\cap N[w]|=1$. Thus,
    \begin{align*}|N[u]\cup N[v]\cup N[w]|&=|N[u]|+|N[v]|+|N[w]|\\ & \quad -|N[u]\cap N[v]| -|N[u]\cap N[w]|-|N[v]\cap N[w]|+|N[u]\cap N[v]\cap N[w]|\\ &= 3(n+1)-3\times 2+1=3n-2.\end{align*}
    Actually, $S=\{\mb 0, e_1+e_2, e_1+e_3\}$.
\item If all $|N[u]\cap N[v]|$, $|N[u]\cap N[w]|$ and $|N[v]\cap N[w]|$ are zero, then $|N[u]\cup N[v]\cup N[w]|=3(n+1)$.
\item If at least one of $|N[u]\cap N[v]|$, $|N[u]\cap N[w]|$ and $|N[v]\cap N[w]|$ is not zero and at least one of them is zero, then
\begin{align*}|N[u]\cup N[v]\cup N[w]|&=|N[u]|+|N[v]|+|N[w]|\\ & \quad -|N[u]\cap N[v]| -|N[u]\cap N[w]|-|N[v]\cap N[w]|+|N[u]\cap N[v]\cap N[w]|\\ &\ge 3(n+1)-2\times 2=3n-1.\end{align*}
\end{enumerate}
Thus, $x_3=3n-5$.

\ms

\nt Let us consider $S=\{u_1, u_2, u_3, u_4\}$, where $u_1, u_2, u_3, u_4$ are distinct. Then
\begin{align*}\left|\bigcup_{l=1}^4 N[u_l]\right|&  =\sum_{l=1}^4 |N[u_l])| -\sum_{1\le j<l\le 4}|N[u_j]\cap N[u_l]|+\sum_{1\le h<j<l\le 4}|N[u_h]\cap N[u_j]\cap N[u_l]|\\ & \quad -|N[u_1]\cap N[u_2]\cap N[u_3]\cap N[u_4]|.\end{align*}
\begin{enumerate}[(1).]
\item If only one of $N[u_j]\cap N[u_l]=\varnothing$, then by Lemma~\ref{lem-3neighbor} the third summand is 1 and fourth summand is 0. Then $\left|\bigcup_{l=1}^4 N[u_l]\right|\ge 4n+4-5\times 2+1=4n-5$.

Actually, $S=\{\mb 0, e_1+e_2, e_1+e_3, e_1+e_2+e_3+e_4\}$.
\item If more than one of $N[u_j]\cap N[u_l]=\varnothing$, then the third and fourth summands are 0. Thus, $\left|\bigcup_{l=1}^4 N[u_l]\right|\ge 4n+4-4\times 2=4n-4$.

\item  If all of $N[u_j]\cap N[u_l]\ne \varnothing$, then \begin{align*}\left|\bigcup_{l=1}^4 N[u_l]\right|&= 4n+4-6\times 2+4-|N[u_1]\cap N[u_2]\cap N[u_3]\cap N[u_4]|\\&\ge 4n-4-1=4n-5.\end{align*}

\end{enumerate}
\nt Therefore, $x_4=4n-9$.

\ms\nt Hence, by Theorem~\ref{thm-neighborhood} we have $\xi\ge 2$ when $n=2$; $\xi\ge 3$ when $n=3$; $\xi\ge 5$ when $n=4$; and $\xi\ge 4n-12$ when $n\ge 5$. Thus, we have $str(Q_2)\ge 6$; $str(Q_3)\ge 11$; $str(Q_4)\ge 21$; and $str(Q_n)\ge 2^n+4n-12$ when $n\ge 5$.
\end{proof}


\nt From the proof of Theorem~\ref{thm-Qn-LB} we have $x_1=n$, $x_2=2n-2$, $x_3=3n-5$, and $x_4=4n-9$ for $Q_n$.

\ms
\nt Suppose $x_{i+1}=|N(S)\setminus S|$ for some vertices subset $S$ with $|S|=i+1$. Let $S=\{u_1, \dots, u_{i+1}\}$ with $|S|=i+1$.

\begin{align*}
x_{i+1}+(i+1)=\left|\bigcup_{l=1}^{i+1} N[u_l]\right|& =|N[u_1]|+\left|\bigcup_{l=2}^{i+1} N[u_l]\right|-\left|N[u_1]\cap \left(\bigcup_{l=2}^{i+1} N[u_l]\right)\right|\\
& =|N[u_1]|+\left|\bigcup_{l=2}^{i+1} N[u_l]\right|-\left|\bigcup_{l=2}^{i+1} (N[u_1]\cap N[u_l])\right|\\&\ge (n+1)+(x_i+i)-2i \tag{since $|N[u_1]\cap N[u_l]|\le 2$}
\end{align*}
So $x_{i+1}\ge n+x_i-2i$. Since $x_4=4n-9$, by induction we will get $x_i\ge in+3-(i-1)i$, where $i\ge 4$.

\ms \nt Let $\eta_i=x_i-i+1$, then $\eta_{i+1}\ge \eta_i +n-2i-1$. So $\eta_i$ is increasing when $i\le (n-1)/2$.

\ms\nt Suppose $n=2m$, $m\ge 2$, then $\eta_{m} \ge \eta_{m-1}+2m -2(m-1)-1=\eta_{m-1}+1$. So $\xi\ge \eta_{m}=x_{m}-m+1\ge m(2m)+3-(m-1)m-m+1=m^2+4$. We have $str(Q_{2m})\ge 2^{2m}+m^2+4$.

\ms\nt Suppose $n=2m-1$, $m\ge 2$, then $\eta_{m}\ge \eta_{m-1} +(2m-1)-2(m-1)-1=\eta_{m-1}$. So $\xi\ge \eta_{m}=x_{m}-m+1\ge m^2-m-4$. We have $str(Q_{2m-1})\ge 2^{2m-1}+m^2-m+4$.

\ms\nt Combining with Theorem~\ref{thm-Qn-LB} we have

\begin{theorem}\label{thm-Qn-LB1} For $n\ge 2$
\begin{enumerate}[$1$.]
\item $str(Q_2)\ge 6$; $str(Q_3)\ge 11$; $str(Q_4)\ge 21$;
\item $str(Q_n)\ge 2^n+4n-12$ for $5\le n\le 9$;
\item $str(Q_{2m})\ge 2^{2m}+m^2+4$ for $m\ge 5$.
\item $str(Q_{2m-1})\ge 2^{2m-1}+m^2-m+4$ for $m\ge 6$;
\end{enumerate}
\end{theorem}

\begin{corollary} $str(Q_2)=6$, $str(Q_3)=11$, $str(Q_4)=21$ and $str(Q_5)=40$.\end{corollary}

\begin{proof} Combining with Theorems~\ref{thm-Qn} and \ref{thm-Qn-LB}, we have $str(Q_2)=6$, $str(Q_3)=11$ and $str(Q_4)=21$.
By considering the following labeling of $Q_5$, we have $str(Q_5)=40$.

\centerline{
\begin{tabular}{c||*{8}{c|}|>{\PBS\centering\hspace{0pt}}m{2.3cm}|}\hline
\diaghead(-3,2){\hskip10mm}%
{$Q_2$}{$Q_3$}& 000 & 100& 110& 010 & 001 & 101 & 111 & 011 & \\\hline \hline
\multirow{2}{1cm}{00} & 00000 & 10000 & 11000 & 01000 & 00100 & 10100 & 11100 & 01100 & \\
 & 1 & 32 & 3 & 31 & 26 & 2 & 30 & 6 & 37\\
\hline
\multirow{2}{1cm}{10} & 00010 & 10010 & 11010 & 01010 & 00110 & 10110 & 11110 & 01110 &\\
& 21 & 4 & 29 & 7 & 12 & 27 & 9 & 24 & 39\\\hline
\multirow{2}{1cm}{11} & 00011 & 10011 & 11011 & 01011 & 00111 & 10111 & 11111 & 01111 &\\
& 16 & 19 & 11 & 20 & 17 & 13 & 18 & 15 & 36\\\hline
\multirow{2}{1cm}{01} & 00001 & 10001 & 11001 & 01001 & 00101 & 10101 & 11101 & 01101 & \\
& 22 & 5 & 28 & 8 & 14 & 25 & 10 & 23 & 39\\\hline\hline
 & 38 & 37 & 40 & 39 & 40 & 40 & 40 & 39 & max. induced edge label\\\hline
 \end{tabular}}

\ms\nt Note that $Q_5\cong Q_3\times Q_2$. The first row and the first column are vertices of $Q_3$ and $Q_2$, respectively. Following is the corresponding figure.

\ms
\centerline{\epsfig{file=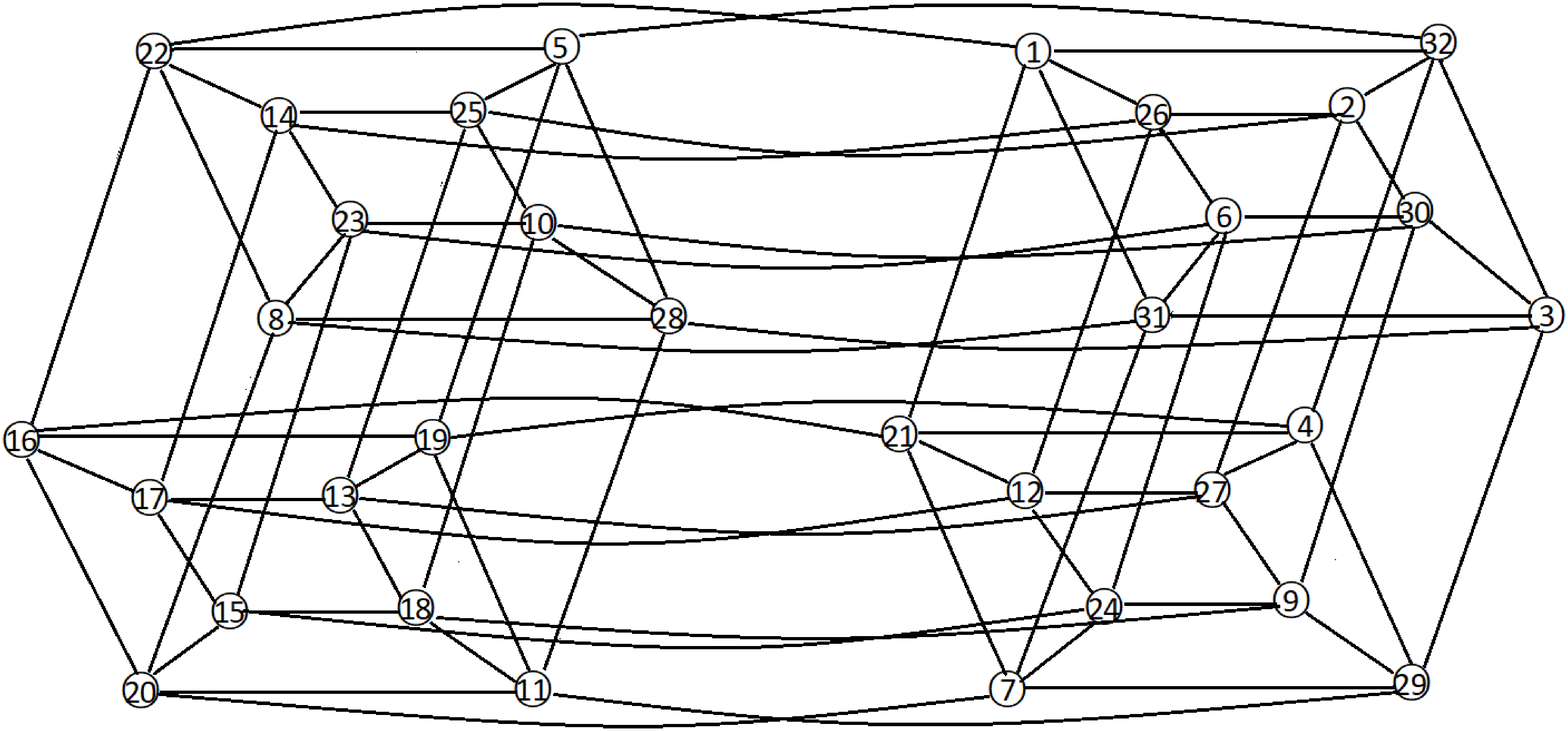, width=15cm}}
\end{proof}

\nt Following is a numbering for $Q_6\cong Q_3\times Q_3$.

\centerline{
\begin{tabular}{c||*{8}{c|}|>{\PBS\centering\hspace{0pt}}m{2.3cm}|}\hline
\diaghead(-3,2){\hskip10mm}%
{$Q_3$}{$Q_3$}& 000 & 100& 110& 010 & 001 & 101 & 111 & 011 & \\\hline \hline
\multirow{2}{1cm}{000} & 000000 & 100000 & 110000 & 010000 & 001000 & 101000 & 111000 & 011000 & \\
 & 1 & 64 & 3 & 63 & 62 & 2 & 52 & 7 & 70\\
\hline
\multirow{2}{1cm}{100} & 000100 & 100100 & 110100 & 010100 & 001100 & 101100 & 111100 & 011100 &\\
& 61 & 4 & 45 & 8 & 11 & 57 & 17 & 58 & 75\\\hline
\multirow{2}{1cm}{110} & 000110 & 100110 & 110110 & 010110 & 001110 & 101110 & 111110 & 011110 &\\
& 15 & 49 & 28 & 38 & 56 & 21 & 33 & 19 & 75\\\hline
\multirow{2}{1cm}{010} & 000010 & 100010 & 110010 & 010010 & 001010 & 101010 & 111010 & 011010 & \\
& 60 & 6 & 50 & 10 & 13 & 54 & 25 & 46 & 79\\\hline
\multirow{2}{1cm}{001} & 000001 & 100001 & 110001 & 010001 & 001001 & 101001 & 111001 & 011001 & \\
 & 59 & 5 & 37 & 9 & 12 & 53 & 23 & 48 & 76\\
\hline
\multirow{2}{1cm}{101} & 000101 & 100101 & 110101 & 010101 & 001101 & 101101 & 111101 & 011101 &\\
& 14 & 44 & 30 & 43 & 55 & 20 & 42 & 18 & 75\\\hline
\multirow{2}{1cm}{111} & 000111 & 100111 & 110111 & 010111 & 001111 & 101111 & 111111 & 011111 &\\
& 39 & 26 & 35 & 31 & 22 & 40 & 32 & 34 & 75\\\hline
\multirow{2}{1cm}{011} & 000011 & 100011 & 110011 & 010011 & 001011 & 101011 & 111011 & 011011 & \\
& 16 & 51 & 27 & 41 & 47 & 24 & 36 & 29 & 78\\\hline\hline
 & 76 & 77 & 78 & 73 & 78 & 78 & 77 & 77 & max. induced edge label\\\hline
 \end{tabular}}

\ms\nt Thus, $76\le str(Q_6)\le 79$.

\section{Conclusion and Open problems}
We obtained a sufficient condition for $str(G)=|V(G)|+\d(G)$. Many open problems are solved and new results are obtained immediately. New lower bound of $str(G)$ in term of $\alpha(G)$ is obtained. Consequently, strength of all 2-regular graphs are determined. An approach of obtaining the lower bound of $str(G)$ in term of the neighborhood size of all possible subsets of $V(G)$ is also obtained. This gives us sharp lower bound of $str(Q_n)$ and answers Problem~\ref{pbm-Qn} partially. The following problems arise naturally.

\begin{problem} Find sufficient and/or necessary conditions such that $str(G) = 2p - 2\alpha(G) + 1$ or $str(G) = p + \xi(G)$.  \end{problem}

\begin{problem} Determine the exact strength of all $r$-regular graphs for $r\ge 3$. \end{problem}

\nt Note that for $G=C_{2n+1}, n\ge 1$, $str(G) = 2n+3 = |V(G)|+\d(G) = 2|V(G)|-2\alpha(G)+1$.

\begin{problem} Characterize all graphs $G$ of order $p$ with (i) $str(G) = p+\d(G) = 2p-2\alpha(G)+1$ or (ii) $str(G) = p+\xi(G) = 2p-2\alpha(G)+1$ or (iii) $str(G) = p + \d(G) = p+\xi(G)$. \end{problem}

\nt Note that every 2-regular graph $\mathcal C_k$ that has $k\ge 2$ odd cycles has $str(\mathcal C_k) = 2|V(\mathcal C_k)| - 2\alpha(\mathcal C_k) + 1 > |V(\mathcal C_k)| + \d(\mathcal C_k)$. Observe that if $\mathcal C_k$ contains an even cycle $C$, then it contains a component $C$ with $str(C) = |V(C)| + \d(C)$.

\begin{problem} Prove that if $str(G) = 2|V(G)| - 2\alpha(G) + 1 > |V(G)| + \d(G)$, then $G$ contains a proper subgraph $H$ with $str(H) = |V(H)| + \d(H)$. \end{problem}

\begin{problem} Prove that for each graph $G$, if $str(G) > |V(G)| + \d(G)$, then $str(G) = 2|V(G)| - 2\alpha(G) + 1$. Otherwise, either $G$ is a proper subgraph of a graph $H$ with $str(H) = 2|V(H)| - 2\alpha(H) + 1$ with $\alpha(H)\ge \alpha(G)$, or else, $G$ contains a proper subgraph $H$ with $str(H) = |V(H)| + \d(H)$. \end{problem}


\begin{thebibliography}{99}

\bibitem{Acharya+Hegde91} B.D. Acharya and S.M. Hegde, Strongly indexable graphs, {\it Discrete Math.}, {\bf 93} (1991), 123--129.

\bibitem{Avadayappan+Jeyanthi+Vasuki01} S. Avadayappan, P. Jeyanthi and R. Vasuki, Super magic strength of a graph, {\it Indian J. Pure Appl. Math.} {bf 32}, (2001) 1621--1630.

\bibitem{Bondy} J.A. Bondy, U.S.R. Murty, {\it Graph theory with applications}, New York, MacMillan, 1976.

\bibitem{Enomoto+Llado+Nakamigawa98} H. Enomoto, A. Llad\'o, T. Nakamigawa, and G. Ringel, Super edge-magic graphs, {\it SUT J. Math.} {\bf 34} (1998), 105--109.

\bibitem{Ichishima+MB+Oshima18} R. Ichishima, F.A. Muntaner-Batle, A. Oshima,  Bounds for the strength of graphs, {\it Aust. J. Combin.} {\bf72(3)},  (2018) 492--508.


\bibitem{Ichishima+MB+Oshima19} R. Ichishima, F.A. Muntaner-Batle, A. Oshima, On the strength of some trees, {\it AKCE Int. J. Graphs Comb.} (Online 2019) doi.org/10.1016/j.akcej.2019.06.002.

\bibitem{Shiu+Wong10} W.C. Shiu and F.S. Wong, Strong vertex-graceful labelings for some double cycles, {\it Congr. Numer.}, {\bf 202}, (2010) 17--24.






\end{thebibliography}
\end{document}